\documentclass[12pt]{article}

\usepackage[all]{xypic}
\usepackage{amssymb}
\usepackage{amsmath}

\usepackage{geometry}
\def\defi{\begin{definition}}
\def\edefi{\end{definition}}
\def\teo{\begin{theorem}}
\def\eteo{\end{theorem}}
\def\prop{\begin{proposition}}
\def\eprop{\end{proposition}}
\def\lema{\begin{lemma}}
\def\elema{\end{lemma}}
\def\coro{\begin{corollary}}
\def\ecoro{\end{corollary}}

\def\downy{\downarrow\kern-.14cm}

\newcommand{\Mnd}{\mathbf{Mnd}}
\newcommand{\Adj}{\mathbf{Adj}}
\newcommand{\Alg}{\mathrm{Alg}}
\newcommand{\Inc}{\mathrm{Inc}}
\newcommand{\I}{\mathrm{I}}
\newcommand{\Aa}{\mathrm{A}}
\newcommand{\F}{\scriptscriptstyle{F}}

\newcommand{\Hh}{\scriptscriptstyle{H}}
\newcommand{\M}{\scriptscriptstyle{M}}
\newcommand{\N}{\scriptscriptstyle{N}}
\newcommand{\Ll}{\scriptscriptstyle{L}}
\newcommand{\R}{\scriptscriptstyle{R}}
\newcommand{\D}{\scriptscriptstyle{D}}
\newcommand{\U}{\scriptscriptstyle{U}}

\newcommand{\K}{\scriptscriptstyle{K}}
\newcommand{\E}{\scriptscriptstyle{E}}
\newcommand{\J}{\scriptscriptstyle{J}}
\newcommand{\f}{\scriptscriptstyle{f}}
\newcommand{\h}{\scriptscriptstyle{h}}

\newcommand{\scirc}{\scriptstyle{\circ}\textstyle{}}
\newcommand{\separ}{\phantom{.}\!}
\newcommand{\mC}{\scriptscriptstyle{\mathcal{C}}}

\def\malg[#1]{\widehat{\textup{Alg}}\textup{-}\mathbb{#1}}
\def\alg[#1]{\mathbb{D}\textup{-Alg}}

\def\Label#1{\label{#1}\ifmmode\llap{[#1] }\else 
\marginpar{\smash{\hbox{\tiny [#1]}}}\fi} 
\def\Label{\label}

\title{Applications of the Kleisli and Eilenberg-Moore 2-adjunctions}

\usepackage{authblk}

\author[2]{L\'opez Hern\'andez, J.L. \thanks{jululohe@gmail.com}}
\author[1]{Turcio Cuevas, L.J. \thanks{luisjturcio@ciencias.unam.mx}}
\author[1]{V\'azquez-M\'arquez, A. \thanks{adrian.vazquez.mqz@gmail.com,$\:$avazquez@uiwbajio.mx}}
\affil[1]{Instituto de Matem\'aticas, UNAM y Universidad Incarnate Word, Campus Baj\'io}
\affil[2]{Facultad de Ciencias, UNAM}

\begin{document}

\maketitle

\abstract{We collect some isomorphisms of categories and bijections of
  structures related to classical monad theory using Kleisli and Eilenberg-Moore 2-adjunctions.}

\section{Introduction}

Motivated by \cite{brto_eicb} and \cite{cljs_klem}, the authors apply
2-adjunctions of Kleisli and Eilenberg-Moore in order to get some classical
isomorphisms of categories and bijections of structures related to monads. \\

Among the examples given in this article, there is one of high importance. In
\cite{moie_motc}, I. Moerdijk gave an equivalence between the lifting of a
monoidal structure, over a category $\mathcal{C}$, to a monoidal structure on
the category of Eilenberg-Moore algebras $\mathcal{C}^{\F}$, for a monad with
endofunctor $F$ on the category $\mathcal{C}$, and the colax monad structures
on $F$ for the monoidal category $\mathcal{C}$. This equivalence of structures
lacks of naturality but using the 2-adjunction of Eilenberg-Moore it can be incorporated.\\

Analogously, the following case is analysed. The equivalence between extensions of monoidal structure over a category
$\mathcal{C}$ to a monoidal structure on the Kleisli category
$\mathcal{C}_{\F}$, for a monad with endofunctor $F$ on the category
$\mathcal{C}$, and the lax monad structures on $F$ for the monoidal category
$\mathcal{C}$, cf. \cite{moie_motc} and \cite{zama_fotm}.\\

The 2-adjunctions of Kleisli and Eilenberg-Moore are generalized to the
context of 2-categories that accept the constructions of algebras. \\

We give the structure of the article.\\

In Section 2, we give the formal 2-adjunction corresponding to the Kleisli situation. In Section 3, we give the formal 2-adjunction corresponding to the Eilenberg-Moore case.\\

In Section 4, we apply the 2-adjunction of EM to the case where the 2-category is ${}_{2}Cat$. In Section 5, we prove the theorem of I. Moerdijk on the
equivalence of lifted monoidal structures and colax monads.\\

In Section 6, we use the Kleisli 2-adjunction for the ${}_{2}Cat$ case. In
Section 7, we apply this 2-adjunction to \emph{extensions} of a monoidal structure on the Kleisli category and relate it with lax monads.\\

In Section 8, we apply the 2-adjunction of Eilenberg-Moore to the very known case of liftings of functors and commutative diagrams for the forgetful functor, check \cite{bofr_haca_II} and \cite{tami_psdl}.\\
 
In Section 9, we relate actions of the category $\mathcal{C}$ over its Kleisli
category $\mathcal{C}_{\F}$ with strong monads.\\

In Section 10, we finalize with left and right functor algebras for a monad and relate this to certain liftings and extensions, respectively, for the underlying functors, cf. \cite{duej_kaee}.\\

We give some remarks on notation. Suppose that we had an adjunction of the
form $\mathfrak{L}\dashv\mathfrak{R}$, then the unit and counit for this
adjunction will be denoted as $\eta^{\mathfrak{R}\mathfrak{L}}$ and
$\varepsilon^{\mathfrak{L}\mathfrak{R}}$, respectively. This notation is complicated but it is clear and prevents the proliferation of several greek letters
to denote new units and counits. As the article develops, the reader might see the advantage in the usage of this notation.\\

We will be working with monoidal categories denoted as $(\mathcal{C}, \otimes, I, a, l, r)$ and also as $(\mathcal{C}, \otimes, I)$ , as a contraction, that leaves understood the natural constrain transformations. We will be working with the constant functor $\delta_{I}:\mathbf{1}\longrightarrow \mathcal{C}$, on $I$, where $\mathbf{1}$ is the category with only one object $0$ and only one arrow $1_{0}$. That is to say, $\delta_{I}(0) = I$.\\

On the other hand, it is known that a category with binary products and a terminal object has a canonical (cartesian) monoidal structure. This is the case for the category $Cat$, of small categories. The natural constraint transformations, taken on components, are functors, for example, for $\mathcal{C}, \mathcal{D}, \mathcal{E}$, $a_{\mathcal{C},\mathcal{D}\!,\mathcal{E}}:(\mathcal{C}\times\mathcal{D})\times\mathcal{E}\longrightarrow \mathcal{C}\times(\mathcal{D}\times\mathcal{E})$ is the obvious functor. In order to compact the notation, we will agree that in the case that the component be the object $\mathcal{C}, \mathcal{C}, \mathcal{C}$, the asociativity functor will be denoted simply as $a_{\mC}$. In turn, the respective constraint functors will be denoted as $l_{\mC}$ and $r_{\mC}$. \\

Finally, the horizontal composition in a general 2-category $\mathcal{A}$ will be denoted as $\cdot$ or by juxtaposition, this notation will be used indistinctively. The vertical composition on 2-cells will be given the symbol $\circ$.

\section{Formal Kleisli 2-Adjunction}

Consider a 2-category such that, its dual, $\mathcal{A}^{op}$ admits the construction of algebras. Due to this property of the 2-category $\mathcal{A}^{op}$, we will be able to construct a 2-adjunction of the form

\begin{equation}\label{1405261543}
\xymatrix@C=1.6cm{
\Mnd(\mathcal{A}^{op}) \ar@<-3pt>[r]_{\Psi_{\K}} & \Adj_{\R}(\mathcal{A}^{op})\ar@<-3pt>[l]_{\Phi_{\K}} 
}
\end{equation}\\

If we describe the 2-adjunction over $\mathcal{A}$ rather than on the  opposite one
then the 2-category $\Mnd(\mathcal{A}^{op})$ will be isomorphic to 
$\Mnd^{\scriptscriptstyle{\bullet}}(\mathcal{A})$ and the 2-category
$\Adj_{\R}(\mathcal{A}^{op})$ will be isomorphic to $\Adj_{\Ll}(\mathcal{A})$. Note that in \cite{stro_fotm}, the category $\mathcal{A}^{op}$ is denoted as $\mathcal{A}^{\ast}$. \\

The description of the 2-category $\Mnd^{\scriptscriptstyle{\bullet}}(\mathcal{A})$ is given as follows.

\begin{enumerate}

\item [1.-] The 0-cells are monads in $\mathcal{A}$, \emph{i.e.} $(A, f,\mu^{\f}, \eta^{\f})$. The short notation $(A,f)$ will be used for such a  monad.  

\item [2.-] The 1-cells, which we call indistinctively as morphisms of  monads, are pairs of the form $(m, \pi):(A, f)\longrightarrow (B, h)$; where  $m:A\longrightarrow B$ is a 1-cell in $\mathcal{A}$ and $\pi:mf\longrightarrow hm$ is a 2-cell in $\mathcal{A}$ such that the following diagrams commute 

\begin{equation*}
\begin{array}{ccc}
\xymatrix@C=1.5cm@R=1.5cm{
mff \ar[r]^-{\pi f}\ar[d]_{m\mu^{\f}} & hmf\ar[r]^-{h \pi } &
hhm\ar[d]^{\mu^{\h} m}\\
mf\ar[rr]_{\pi} & &     hm
} & 
\xy<1cm,0cm>
\POS (0,-20.5) *+{,},
\endxy
&
\xy<1cm,0cm>
\POS (0, -6)
\xymatrix@C=1cm@R=0.8cm{
& m\ar[ld]_{m\eta^{f}}\ar[rd]^{\eta^{h} m}& \\
mf\ar[rr]_{\pi} & & hm
}
\endxy
\end{array}
\end{equation*}

\item [3.-] The 2-cells, which we call indistinctively as transformations of  monads, have the form $\vartheta:(m, \pi)\longrightarrow (n, \tau):(A, f)\longrightarrow (B, h)$, such that $\vartheta:m\longrightarrow n:A\longrightarrow B$ is a 2-cell in $\mathcal{A}$ and the following diagram commutes

\begin{equation*}\Label{1405261627}
\xy<1cm,0cm>
\POS (0, 0) *+{mf} = "c1",
\POS (20, 0) *+{hm} = "c2", 
\POS (0, -20) *+{nf} = "d1",
\POS (20, -20) *+{hn} = "d2", 
\POS (18, -10) = "a",
\POS (10, -18) = "b",
\POS "c1" \ar^{\pi} "c2",
\POS "d1" \ar_{\tau} "d2",
\POS "c1" \ar_{\vartheta\! f} "d1",
\POS "c2" \ar^{h \vartheta } "d2",
\endxy
\end{equation*}

This 2-cell is displayed as follows 

\begin{equation*}
\xy<1cm,0cm>
\POS (0, 0) *+{(A, f)} = "c1",
\POS (32, 0) *+{(B, h)} = "c2", 
\POS (16,0) *+{\vartheta},
\POS (10, 3.5) = "a",
\POS (10, -3.5) = "b",
\POS "c1" \ar@/^1.5pc/^{(m,\: \pi)} "c2",
\POS "c1" \ar@/_1.5pc/_{(n,\:\tau)} "c2",
\POS "a" \ar "b",
\endxy
\end{equation*}

\end{enumerate}

The structure of the 2-category $\Adj_{\Ll}(\mathcal{A})$ is given as follows

\begin{enumerate}

\item [1.-] The $0$-cells are made of adjunctions

\begin{equation*}\Label{1407051035}
\begin{array}{cc}
\xymatrix@C=1.6cm{
A \ar@<-3pt>[r]_{l} & B\ar@<-3pt>[l]_{r} 
} &\!\!\! .
\end{array}
\end{equation*}

\item [2.-] The $1$-cells are of the form $(j, k, \rho)$ such that the second diagram is the 2-cell mate of the first one that commutes

\begin{equation*}\Label{1407051036}
\begin{array}{ccc}
\xy<1cm,0cm>
\POS (0, 0) *+{A} = "c1",
\POS (20, 0) *+{\overline{A}} = "c2", 
\POS (0, -20) *+{B} = "d1",
\POS (20, -20) *+{\overline{B}} = "d2", 
\POS (18, -10) = "a",
\POS (10, -18) = "b",
\POS "c1" \ar^{j} "c2",
\POS "c1" \ar_{l} "d1",
\POS "c2" \ar^{\overline{l}} "d2",
\POS "d1" \ar_{k} "d2",
\endxy & 
\xy<1cm,0cm>
\POS (0,-22) *+{,},
\endxy
& 
\xy<1cm,0cm>
\POS (0, 0) *+{A} = "c1",
\POS (20, 0) *+{\overline{A}} = "c2", 
\POS (0, -20) *+{B} = "d1",
\POS (20, -20) *+{\overline{B}} = "d2",
\POS (10, -2) = "b",
\POS (18, -10) = "a",
\POS "c1" \ar^{j} "c2",
\POS "d1" \ar^{r} "c1",
\POS "d2" \ar_{\overline{r}} "c2",
\POS "d1" \ar_{k} "d2",
\POS "b" \ar@/_.5pc/_{\rho}"a",
\endxy 
\end{array}
\end{equation*}

The mate $\rho$ is described, since the left one commutes, by 

\begin{equation}\Label{1407052237}
\rho = \overline{r}k\varepsilon\circ\overline{\eta}jr
\end{equation}

This morphism can be represented as

\begin{equation*}\Label{1407052238}
\xy<1cm,0cm>
\POS (0, 0) *+{A} = "c1",
\POS (20, 0) *+{\overline{A}} = "c2", 
\POS (0, -20) *+{B} = "d1",
\POS (20, -20) *+{\overline{B}} = "d2",
\POS "c1" \ar^{j} "c2",
\POS "c1" \ar@<-2pt>_{l} "d1",
\POS "d1" \ar@<-2pt>_{r} "c1",
\POS "c2" \ar@<-2pt>_{\overline{l}} "d2",
\POS "d2" \ar@<-2pt>_{\overline{r}} "c2",
\POS "d1" \ar_{k} "d2",
\POS (10, -2) = "b",
\POS (16, -9) = "a",
\POS "b" \ar@/_.4pc/_{\rho}"a",
\endxy
\end{equation*}

\noindent and denoted as  $(j,k,\rho):l\dashv r\longrightarrow
\overline{l}\dashv \overline{r}$. Since the diagram corresponding to the left
adjoints commutes, the 2-category of adjunctions has the subindex $L$.

\item [3.-] The $2$-cells are made of a pair of 2-cells in $\mathcal{A}$, $(\alpha, \beta)$ as in  

\begin{equation*}\Label{1405141814}
\xy<1cm,0cm>
\POS (0, 0) *+{A} = "c1",
\POS (25, 0) *+{\overline{A}} = "c2", 
\POS (0, -25) *+{B} = "d1",
\POS (25, -25) *+{\overline{B}} = "d2", 
\POS (12.5, 0) *+{\alpha} = "e1",
\POS (12.5, -25) *+{\beta} = "e2",
\POS (7, 3) = "a",
\POS (7, -2.5) = "b",
\POS (7,-22) = "ap",
\POS (7, -27.5) = "bp",
\POS "c1" \ar@/^1.2pc/^{j} "c2",
\POS "c1" \ar@/_1.2pc/_{j'} "c2",
\POS "c1" \ar@<-2pt>_{l} "d1",
\POS "d1" \ar@<-2pt>_{r} "c1",
\POS "c2" \ar@<-2pt>_{\overline{l}} "d2",
\POS "d2" \ar@<-2pt>_{\overline{r}} "c2",
\POS "d1" \ar@/^1.2pc/^{k} "d2",
\POS "d1" \ar@/_1.2pc/_{k'} "d2",
\POS "a" \ar "b",
\POS "ap" \ar "bp",
\endxy
\end{equation*}

\noindent such that they fulfill one of the following equivalent conditions 

\begin{enumerate}

\item [(i)] $\overline{l}\alpha = \beta l$,

\item  [(ii)]  $\rho'\circ\alpha r = \overline{r}\beta\circ\rho$.\\

\end{enumerate}

\newtheorem{1407060859}{Remark}[section]
\begin{1407060859}
Note that the previous conditions can be seen as commutative surface diagrams.
\end{1407060859}

This $2$-cell can be displayed as follows

\begin{equation*}
\xy<1cm,0cm>
\POS (0, 0) *+{l\dashv r} = "c1",
\POS (32, 0) *+{\overline{l}\! \dashv \! \overline{r}} = "c2", 
\POS (16,0) *+{\scriptstyle{(}\scriptstyle{\alpha},\scriptstyle{\beta})},
\POS (10, 3.5) = "a",
\POS (10, -3.5) = "b",
\POS "c1" \ar@/^1.5pc/^{(j, \:k,\:\rho)} "c2",
\POS "c1" \ar@/_1.5pc/_{(j', \:k'\!,\:\rho')} "c2",
\POS "a" \ar "b",
\endxy
\end{equation*}

The $n$-cell structure described arrange itself to form a $2$-category, that
is to say inherits the $2$-category structure of $\mathcal{A}$. 

\end{enumerate}

Before going into the details on the construction of the 2-functor $\Psi_{\K}$, we develop some
calculations. These calculations are dual to those made in
\cite{stro_fotm}. Note that we are going to be switching between the
2-categories $\mathcal{A}^{op}$ and $\Mnd(\mathcal{A}^{op})$ to $\mathcal{A}$
and $\Mnd^{\scriptscriptstyle{\bullet}}(\mathcal{A})$, respectively. \\

Since the 2-category $\mathcal{A}^{op}$ admits the construction of algebras,
the functor $\Inc_{\mathcal{A}^{op}}:\mathcal{A}^{op}\longrightarrow \Mnd(\mathcal{A}^{op})$ admits a right adjoint, denoted as
$\Alg_{\mathcal{A}^{op}}:\Mnd(\mathcal{A}^{op})\longrightarrow
\mathcal{A}^{op}$. These 2-functors are going to be short denoted as
$\I^{\scriptscriptstyle{\bullet}}$ and $\Aa^{\scriptscriptstyle{\bullet}}$,
respectively. \\

The corresponding counit, on the component $(A, f^{op})$, is $\varepsilon^{I\!
  A^{\scriptscriptstyle{\bullet}}}(A,
f^{op}):\Inc_{\mathcal{A}^{op}}\Alg_{\mathcal{A}^{op}}(A,
f^{op})\longrightarrow (A, f^{op})$. If we define $\Alg_{\mathcal{A}^{op}}(A,
f^{op}) = A_{\f}$, the \emph{Kleisli object}, then $\varepsilon^{I\!  A^{\scriptscriptstyle{\bullet}}}(A,f^{op}) = (g_{\f}, \iota_{\f}): (A, f)\longrightarrow (A_{\f},1_{A_{\f}})$. This last 1-cell belongs to $\Mnd^{\scriptscriptstyle{\bullet}}(\mathcal{A})$, where $g_{\f}:A\longrightarrow A_{\f}$ and $\iota_{\f}:g_{\f}f\longrightarrow g_{\f}$.\\

Following \cite{stro_fotm}, for any monad $(A, f^{op})$ in $\Mnd(\mathcal{A}^{op})$, there exists an adjunction in $\mathcal{A}$,

\begin{equation*}\Label{1406142138}
\begin{array}{cc}
\xymatrix@C=1.6cm{
A \ar@<-3pt>[r]_{g_{\f}} & A_{\f}\ar@<-3pt>[l]_{v_{\f}} 
} &\!\!\! ,
\end{array}
\end{equation*}\\

\noindent such that it generates the monad $(A, f)$, with unit $\eta^{\f}$ and
counit $\varepsilon^{\separ gv_{\!\f}}$. It can be checked that $\iota_{\f} =
\varepsilon^{\separ gv_{\!\f}}g_{\f}$. This adjunction is called the \emph{Kleisli adjunction}.\\

Suppose that there is a morphism of monads $(m^{op}, \pi):(B, h^{op})\longrightarrow (A, f^{op})$ in $\Mnd(\mathcal{A}^{op})$, \emph{i.e.} $(m,\pi):(A,f)\longrightarrow (B,h)$ in $\Mnd^{\scriptscriptstyle{\bullet}}(\mathcal{A})$. Take the following composition of morphisms of monads $(g_{\h},\iota_{\h})\cdot(m,\pi)= (g_{\h}m, \iota_{\h}m\circ g_{\h}\pi):(A,f)\longrightarrow (B_{\h}, 1_{B_{\h}})$.\\

Since the counit is universal from $\Inc_{\mathcal{A}^{op}}$ to $(A,f^{op})$, there exists a 1-cell $m_{\pi}:A_{\f}\longrightarrow B_{\h}$, in $\mathcal{A}$, such that the following diagram commute\\

\begin{equation*}\Label{1406021518}
\xy<1cm,0cm>
\POS (25,15) *+{(A, f)} = "a12"
\POS (0,0) *+{(A_{\f}, 1_{A_{\f}})} = "a21",
\POS (50,0) *+{(B_{\h}, 1_{B_{\h}})} = "a23"

\POS "a12" \ar_{(g_{\f},\iota_{\f})} "a21",
\POS "a12" \ar^{\phantom{12}(g_{\h}m, \separ\iota_{\h}m\separ\circ\separ g_{\h}\pi)} "a23",
\POS "a21" \ar_{(m_{\pi}, 1_{m_{\pi}})} "a23",
\endxy
\end{equation*}\\

In particular, $g_{\h}m = m_{\pi}g_{\f}$ and $\iota_{\h}m\circ g_{\h}\pi = m_{\pi}\iota_{\f}$. Note that the associated mate to the first equality is $\rho_{\pi} = v_{\h}m_{\pi}\varepsilon^{g v_{\h}}\circ\eta^{\h}m v_{\f}$ and that $\rho_{\pi} g_{\f} = \pi$.\\

Consider a 2-cell of monads $\vartheta: (m,\pi)\longrightarrow (n,\tau):(A,f)\longrightarrow (B,h)$ in $\Mnd^{\scriptscriptstyle{\bullet}}(\mathcal{A})$. Due to the construction of algebras for $\mathcal{A}^{op}$, the 2-adjunction $\Alg_{\mathcal{A}^{op}}\dashv \Inc_{\mathcal{A}^{op}}$ provides an isomorphism of categories, for $(A,f^{op})$ in $\Mnd(\mathcal{A}^{op})$ and $B$ in $\mathcal{A}^{op}$, of
the form

\begin{equation*}\Label{1406191447}
Hom_{\Mnd(\mathcal{A}^{op})}\big((A,f^{op}), \Inc_{\mathcal{A}^{op}}(B)\big)
\cong Hom_{\mathcal{A}^{op}}\big(\Alg_{\mathcal{A}^{op}}(A,f^{op}), B\big)
\end{equation*}\\

\noindent this translates, in the non-opposite case, into the following assignment 

\begin{equation}\Label{1406191526}
\begin{array}{ccc}
\xy<1cm,0cm>
\POS (0,0) *+{A_{\f}} = "a11",
\POS (25,0) *+{B} = "a12",
\POS (9, 3.5) = "a",
\POS (9,-3.5) = "b"

\POS "a11" \ar@/^1.3pc/^{a} "a12",
\POS "a11" \ar@/_1.3pc/_{b} "a12",

\POS "a" \ar^{\phantom{1}\alpha} "b",

\endxy & \longmapsto &
\xy<1cm,0cm>
\POS (0,0) *+{(A,f)} = "a11",
\POS (30,0) *+{(B, 1_{B})} = "a13",
\POS (9, 3.5) = "a",
\POS (9,-3.5) = "b",

\POS (15, -0.4) *+{\scriptstyle{\alpha\separ g_{\!\f}}} = "a12",

\POS "a11" \ar@/^1.5pc/^{(ag_{\!\f},\separ a\iota_{\!\f})} "a13",
\POS "a11" \ar@/_1.5pc/_{(b\separ g_{\!\f}, \separ b\separ\iota_{\!\f})} "a13",

\POS "a" \ar "b",

\endxy
\end{array}
\end{equation}\\

\noindent On the other hand, we have an equality of
 2-cells\\

\begin{equation*}\Label{1406191527}
\begin{array}{ccc}
\xy<1cm,0cm>
\POS (0,0) *+{(A,f)} = "a11",
\POS (30,0) *+{\phantom{123}(B_{\h}, 1_{B_{\h}})} = "a13",
\POS (10, 3.5) = "a",
\POS (10,-3.5) = "b",

\POS (15, 0.5) *+{\scriptstyle{g_{\h}\vartheta}} = "a12",

\POS "a11" \ar@/^1.5pc/^{(g_{\h}m, \separ \iota_{\h}m\circ g_{\h}\pi)} "a13",
\POS "a11" \ar@/_1.5pc/_{(g_{\h}n, \separ \iota_{\h} n\circ g_{\h}\tau)} "a13",

\POS "a" \ar "b",

\endxy & = &
\xy<1cm,0cm>

\POS (0,0) *+{(A, f)} = "a11",
\POS (30,0) *+{\phantom{123}(B_{\h}, 1_{B_{\h}})} = "a13",
\POS (11, 3.5) = "a",
\POS (11,-3.5) = "b",

\POS (16, 0.5) *+{\scriptstyle{g_{\h}\vartheta}} = "a12",

\POS "a11" \ar@/^1.5pc/^{(m_{\pi}g_{\!\f} , \separ m_{\pi}\iota_{\!\f})} "a13",
\POS "a11" \ar@/_1.5pc/_{(n_{\tau}g_{\!\f} , \separ n_{\tau}\iota_{\!\f})} "a13",

\POS "a" \ar "b",

\endxy
\end{array}
\end{equation*}\\

Therefore, to the 2-cell $g_{\h}\vartheta$ there corresponds, through the
asignment \eqref{1406191526}, a 2-cell $\beta_{\vartheta} =
\Alg_{\mathcal{A}^{op}}(g_{\h}\vartheta)\cdot
\eta^{IA^{\scriptscriptstyle{\bullet}}}(B_{\h})$, such that $g_{\h}\vartheta =
\beta_{\vartheta}g_{\f}$, where $\beta_{\vartheta}:m_{\pi}\rightarrow n_{\tau}$. We
change, at this point, the notation as $\beta_{\vartheta} = \widetilde{\vartheta}$.\\

Without any further ado, we provide the description of the 2-functor $\Psi_{\K}$ \\

\begin{enumerate}

\item [1.-] For the monad $(A, f, \mu^{\f}, \eta^{\f})$ in
  $\Mnd^{\scriptscriptstyle{\bullet}}(\mathcal{A})$, $\Psi_{\K}(A, f) =
  g_{\f}\dashv v_{\f}$, \emph{i.e.} the Kleisli adjunction.

\item [2.-] For the morphism $(m, \pi):(A,f)\longrightarrow (B, h)$,  $\Psi_{\K}(m, \pi) = (m, m_{\pi}, \rho_{\pi})$

\item [3.-] For the transformation $\vartheta: (m,\pi)\longrightarrow (n,\tau):(A,f)\longrightarrow(B,g)$, $\Psi_{\K}(\vartheta) =(\vartheta, \widetilde{\vartheta})$, where $\widetilde{\vartheta}$ is given as above.

\end{enumerate}

\noindent The description of the functor $\Phi_{\K}$ is given as follows

\begin{enumerate}

\item [1.-] For the adjunction $l\dashv r$, $\Phi_{\K}(l\dashv r) = (A, rl)$

\item [2.-] For the morphism of adjunctions $(j,k,\rho):(l\dashv r)\longrightarrow(\overline{l}\dashv \overline{r})$, $\Phi_{\K}(j,k,\rho) =  (j,\pi_{\rho})$. Where $\pi_{\rho} = \rho\separ l$.

\item [3.-] For the transformation of adjunctions $(\alpha, \beta):(j,k,\rho)\longrightarrow (j',k',\rho'):l\dashv r\longrightarrow \overline{l}\dashv\overline{r}$, $\Phi_{\K}(\alpha, \beta) = \vartheta_{(\alpha,\beta)} = \alpha$.

\end{enumerate}

Yet again, following \cite{stro_fotm}, it can be shown that for the adjunction $l\dashv r$, there exists a \emph{dual comparison 1-cell} $k_{rl}:A_{rl}\longrightarrow B$, such that $l = k_{rl}g_{rl}$, $v_{rl} = rk_{rl}$ and $\varepsilon^{rl} l = k_{rl}\iota_{rl}$.\\

The unit of the 2-adjunction in \eqref{1405261543}, $\eta^{\Phi\Psi_{\K}}:1_{\Mnd^{\scriptscriptstyle{\bullet}}(\mathcal{A})}\longrightarrow\Phi_{\K}\Psi_{\K}$ is defined, in the component $(A, f)$, as follows

\begin{equation*}\Label{1406281141}
\eta^{\Phi\Psi_{\K}}(A, f) := (1_{A}, 1_{\f}):(A, f)\longrightarrow (A,f)
\quad in \quad \Mnd^{\scriptscriptstyle{\bullet}}(\mathcal{A})
\end{equation*}\\

In turn, the counit $\varepsilon^{\Psi\Phi_{\K}}:\Psi_{\K}\Phi_{\K}\longrightarrow 1_{\Adj_{l}(\mathcal{A})}$ is defined, in the component $l\dashv r$, as follows

\begin{equation*}\Label{1406282009}
\varepsilon^{\Psi\Phi_{\K}}(l\dashv r) := (1_{A}, k_{rl}, 1_{v_{rl}}):g_{rl}\dashv v_{rl}\longrightarrow l\dashv r
\quad in \quad \Adj_{\Ll}(\mathcal{A})
\end{equation*}\\

\newtheorem{1406031238}{Theorem}[section]
\begin{1406031238}
There exists a 2-adjunction $\Psi_{\K}\dashv \Phi_{\K}$.
\end{1406031238}

\noindent \emph{Proof}:\\

We prove only one of the triangular identities, \emph{i.e.} $\Phi_{\K}\separ\varepsilon^{\Psi\Phi_{\K}}\circ\eta^{\Phi\Psi_{\K}}\Phi_{\K} = 1_{\Phi_{\K}}$.\\

\begin{eqnarray*}
\big(\Phi_{\K}\separ\varepsilon^{\Psi\Phi_{\K}}\circ\eta^{\Phi\Psi_{\K}}\Phi_{\K}\big)(l\dashv r) &=& \Phi_{\K}\separ\varepsilon^{\Psi\Phi_{\K}}(l\dashv r)\cdot\eta^{\Phi\Psi_{\K}}\Phi_{\K}(l\dashv r) \\
&=& \Phi_{\K}(1_{A}, k_{rl}, 1_{v_{rl}}) \cdot \eta^{\Phi\Psi_{\K}}(A, rl)\\
&=& (1_{A}, 1_{v_{rl}}g_{rl})\cdot (1_{A}, 1_{rl}) = (1_{A}, 1_{rl}) = 1_{(A, rl)}\\
&=& 1_{\Phi_{\K}(l\dashv\separ r)} = 1_{\Phi_{\K}}(l\dashv r).
\end{eqnarray*}

\begin{flushright}
$\square$
\end{flushright}

Since the left 2-adjoint $\Psi_{\K}$ assigns the Kleisli adjunction to a
monad, the 2-adjunction is called \emph{Kleisli 2-adjunction}.


\section{Formal Eilenberg-Moore 2-Adjunction}

Consider a $2$-category $\mathcal{A}$ which admits the construction of algebras. With this property of $\mathcal{A}$, we will construct a 2-adjunction of the form

\begin{equation*}\Label{1405151703}
\begin{array}{cc}
\xymatrix@C=1.6cm{
\Adj_{\R}(\mathcal{A}) \ar@<-3pt>[r]_{\Phi_{\E}} & \Mnd(\mathcal{A})\ar@<-3pt>[l]_{\Psi_{\!\E}} 
} &\!\!\! ,
\end{array}
\end{equation*}\\

\noindent The 2-category $\Adj_{\R}(\mathcal{A})$ is described as follows 

\begin{enumerate}

\item [1.-] The $0$-cells are made of adjunctions

\begin{equation*}\Label{1404151713}
\begin{array}{cc}
\xymatrix@C=1.6cm{
A \ar@<-3pt>[r]_{l} & B\ar@<-3pt>[l]_{r} 
} &\!\!\! .
\end{array}
\end{equation*}

\item [2.-] The $1$-cells are pairs, of 1-cells in $\mathcal{A}$, $(j, k)$
  such that the first diagram is the 2-cell mate to the second commutative one

\begin{equation*}\Label{1405141753}
\begin{array}{ccc}
\xy<1cm,0cm>
\POS (0, 0) *+{A} = "c1",
\POS (20, 0) *+{\overline{A}} = "c2", 
\POS (0, -20) *+{B} = "d1",
\POS (20, -20) *+{\overline{B}} = "d2", 
\POS (18, -10) = "a",
\POS (10, -18) = "b",
\POS "c1" \ar^{j} "c2",
\POS "c1" \ar_{l} "d1",
\POS "c2" \ar^{\overline{l}} "d2",
\POS "d1" \ar_{k} "d2",
\POS "a" \ar@/_.5pc/_{\lambda}"b",
\endxy & 
\xy<1cm,0cm>
\POS (0,-22) *+{,},
\endxy
& 
\xy<1cm,0cm>
\POS (0, 0) *+{A} = "c1",
\POS (20, 0) *+{\overline{A}} = "c2", 
\POS (0, -20) *+{B} = "d1",
\POS (20, -20) *+{\overline{B}} = "d2", 
\POS "c1" \ar^{j} "c2",
\POS "d1" \ar^{r} "c1",
\POS "d2" \ar_{\overline{r}} "c2",
\POS "d1" \ar_{k} "d2",
\endxy 
\end{array}
\end{equation*}

The mate is described by 

\begin{equation}\Label{1405151810}
\lambda = \overline{\varepsilon}kl\circ \overline{l}j\eta
\end{equation}

This morphism can be represented as

\begin{equation*}\Label{1403151811}
\xy<1cm,0cm>
\POS (0, 0) *+{A} = "c1",
\POS (20, 0) *+{\overline{A}} = "c2", 
\POS (0, -20) *+{B} = "d1",
\POS (20, -20) *+{\overline{B}} = "d2", 
\POS (17, -13) = "a",
\POS (11, -19) = "b",
\POS (17,-13) = "ap",
\POS (11, -19) = "bp",
\POS "c1" \ar^{j} "c2",
\POS "c1" \ar@<-2pt>_{l} "d1",
\POS "d1" \ar@<-2pt>_{r} "c1",
\POS "c2" \ar@<-2pt>_{\overline{l}} "d2",
\POS "d2" \ar@<-2pt>_{\overline{r}} "c2",
\POS "d1" \ar_{k} "d2",
\POS "ap" \ar@/_.25pc/_{\lambda}"bp",
\endxy
\end{equation*}

\noindent and denoted as  $(j,k,\lambda):l\dashv r\longrightarrow \overline{l}\dashv \overline{r}$.

\item [3.-] The $2$-cells are made of a pair of 2-cells in $\mathcal{A}$, $(\alpha, \beta)$ as in  

\begin{equation*}\Label{1407190919}
\xy<1cm,0cm>
\POS (0, 0) *+{A} = "c1",
\POS (25, 0) *+{\overline{A}} = "c2", 
\POS (0, -25) *+{B} = "d1",
\POS (25, -25) *+{\overline{B}} = "d2", 
\POS (12.5, 0) *+{\alpha} = "e1",
\POS (12.5, -25) *+{\beta} = "e2",
\POS (7, 3) = "a",
\POS (7, -2.5) = "b",
\POS (7,-22) = "ap",
\POS (7, -27.5) = "bp",
\POS "c1" \ar@/^1.2pc/^{j} "c2",
\POS "c1" \ar@/_1.2pc/_{j'} "c2",
\POS "c1" \ar@<-2pt>_{l} "d1",
\POS "d1" \ar@<-2pt>_{r} "c1",
\POS "c2" \ar@<-2pt>_{\overline{l}} "d2",
\POS "d2" \ar@<-2pt>_{\overline{r}} "c2",
\POS "d1" \ar@/^1.2pc/^{k} "d2",
\POS "d1" \ar@/_1.2pc/_{k'} "d2",
\POS "a" \ar "b",
\POS "ap" \ar "bp",
\endxy
\end{equation*}

\noindent such that they fulfill one of the following equivalent conditions 

\begin{enumerate}

\item [(i)] $\lambda'\circ \overline{l} \alpha = \beta l \circ\lambda$,

\item  [(ii)]  $\alpha r = \overline{r} \beta$.

\end{enumerate}

\newtheorem{1407152212}{Remark}[section]
\begin{1407152212}
Note that the previous conditions can be seen as commutative surface diagrams.
\end{1407152212}

This $2$-cell can be displayed as follows

\begin{equation*}
\xy<1cm,0cm>
\POS (0, 0) *+{l\dashv r} = "c1",
\POS (32, 0) *+{\overline{l}\! \dashv \! \overline{r}} = "c2", 
\POS (16,0) *+{\scriptstyle{(}\scriptstyle{\alpha},\scriptstyle{\beta})},
\POS (10, 3.5) = "a",
\POS (10, -3.5) = "b",
\POS "c1" \ar@/^1.5pc/^{(j, k, \lambda)} "c2",
\POS "c1" \ar@/_1.5pc/_{(j',\separ k',\separ\lambda')} "c2",
\POS "a" \ar "b",
\endxy
\end{equation*}

The described $n$-cell structure arrange itself to form a $2$-category. 

\end{enumerate}

\noindent The 2-category $\Mnd(\mathcal{A})$ is formed as follows

\begin{enumerate}

\item [1.-] The 0-cells are monads in $\mathcal{A}$, $(A,f,\mu^{\f},\eta^{\f})$. The short notation $(A, f)$ will be used for such a monad.

\item [2.-] The 1-cells are \emph{morphisms of monads} $(p,\varphi):(A, f)\longrightarrow (B, h)$ which consist of a 1-cell $p:A\longrightarrow  B$ and a 2-cell $\varphi:hp\longrightarrow pf$, both in $\mathcal{A}$, such that the following  diagrams commutes

\begin{equation*}\Label{1405151917}
\begin{array}{ccc}
\xymatrix@C=1.5cm@R=1.5cm{
hhp \ar[r]^-{h\varphi }\ar[d]_{\mu^{h}\! p} & hpf\ar[r]^-{\varphi f} & pff\ar[d]^{p\mu^{f}}\\
hp\ar[rr]_{\varphi} & &     pf
} & 
\xy<1cm,0cm>
\POS (0,-20.5) *+{,},
\endxy
&
\xy<1cm,0cm>
\POS (0, -7)
\xymatrix@C=1cm@R=0.8cm{
& p\ar[ld]_{\eta^{h}p}\ar[rd]^{p\eta^{f}}& \\
hp\ar[rr]_{\varphi} & & pf
}
\endxy
\end{array}
\end{equation*}

\item [3.-] The 2-cells or \emph{transformations of monads} $\theta:(p,\varphi)\longrightarrow(q,\psi):(A,f)\longrightarrow(B,h)$, consists of a 2-cell $\theta:p\longrightarrow q$ in $\mathcal{A}$ and fulfills  the commutativity of the following diagram

\begin{equation*}\Label{1405151923}
\xy<1cm,0cm>
\POS (0, 0) *+{hp} = "c1",
\POS (20, 0) *+{pf} = "c2", 
\POS (0, -20) *+{hq} = "d1",
\POS (20, -20) *+{qf} = "d2", 
\POS (18, -10) = "a",
\POS (10, -18) = "b",
\POS "c1" \ar^{\varphi} "c2",
\POS "d1" \ar_{\psi} "d2",
\POS "c1" \ar_{h\theta } "d1",
\POS "c2" \ar^{\theta f} "d2",
\endxy
\end{equation*}

This 2-cell is displayed as follows 

\begin{equation*}
\xy<1cm,0cm>
\POS (0, 0) *+{(A, f)} = "c1",
\POS (32, 0) *+{B, h)} = "c2", 
\POS (16,0) *+{\theta},
\POS (10, 3.5) = "a",
\POS (10, -3.5) = "b",
\POS "c1" \ar@/^1.5pc/^{(p,\: \varphi)} "c2",
\POS "c1" \ar@/_1.5pc/_{(q,\:\psi)} "c2",
\POS "a" \ar "b",
\endxy
\end{equation*}

\end{enumerate}

The description of the $2$-functor $\Phi_{\E}$ is given as follows

\begin{enumerate}

\item [1.-] On $0$-cells, $\Phi_{\E}(l\dashv r) = (A, rl, r\varepsilon l, \eta)$, \emph{i.e.} the induced monad by an adjunction.

\item [2.-] On $1$-cells, $(j,k,\lambda):(A, rl)\longrightarrow (\overline{A}, \overline{r}\overline{l})$, $\Phi_{\E}(j,k,\lambda) = (j, \overline{r}\lambda): (A, rl) \longrightarrow (\overline{A}, \overline{r}\overline{l})$.

\item [3.-] On $2$-cells, $(\alpha,\beta):(j,k,\lambda)\longrightarrow(j',k',\lambda')$, $\Phi_{\E}(\alpha,\beta)= \alpha: (j, \overline{r}\lambda)\longrightarrow (j', \overline{r}\lambda')$.

\end{enumerate}

\noindent Before the description of the 2-functor $\Psi_{\E}$, we realize some calculations.\\

Since the 2-category $\mathcal{A}$ admits the construction of algebras, the 2-functor $\Inc_{\mathcal{A}}:\mathcal{A}\longrightarrow \Mnd(\mathcal{A})$ admits a right adjoint, denoted as $\Alg_{\mathcal{A}}:\Mnd(\mathcal{A})\longrightarrow\mathcal{A}$.\\

The corresponding counit, on the component $(A, f)$, is
$\varepsilon^{IA}(A,f):\Inc_{\mathcal{A}}\Alg_{\mathcal{A}}(A,
f)\longrightarrow (A, f)$. If we define $\Alg_{\mathcal{A}}(A,f) = A^{\f}$,
the \emph{Eilenberg-Moore object} for $(A,f)$, then $\varepsilon^{IA}(A,f) :=(u^{\f},\chi^{\f}):(A^{\f}, 1_{A^{\f}})\longrightarrow (A,f)$, where $u^{\f}:A^{\f}\longrightarrow A$ and $\chi^{\f}:u^{\f}d^{\f}u^{\f}\longrightarrow u^{\f}$.\\

In Theorem 2, at \cite{stro_fotm}, the author proved that if $\mathcal{A}$ admits the construction of algebras then for any monad $(A, f)$ in $\Mnd(\mathcal{A})$, there exists an adjunction in $\mathcal{A}$

\begin{equation*}\Label{1405151952}
\begin{array}{cc}
\xymatrix@C=1.6cm{
A \ar@<-3pt>[r]_{d^{\f}} & A^{\f}\ar@<-3pt>[l]_{u^{\f}} 
} &\!\!\! ,
\end{array}
\end{equation*}

\noindent such that it generates the monad $(A, f)$, with unit $\eta^{\f}$ and
counit $\varepsilon^{d\separ u^{\f}}$. It can be checked that $\chi^{\f} =
u^{\f}\varepsilon^{d\separ u^{\f}}$. This adjunction is called the \emph{The
  Eilenberg-Moore adjunction}.\\

Suppose there is a morphism of monads $(p, \varphi):(A, f)\longrightarrow (B,h)$. Take the composition of morphisms of monads $(p,\varphi)\cdot(u^{\f},\chi^{\f})= (pu^{\f}, p\chi^{\f}\circ\varphi u^{\f}):\Inc_{\mathcal{A}}(A^{\f}) = (A^{\f},1_{A^{\f}})\longrightarrow (B, h)$.\\

The previous counit, $\varepsilon^{I\! A}$, is universal from the functor $\Inc_{\mathcal{A}}$, in particular, for the 1-cell 
$(pu^{\f}, p\chi^{\f}\circ\varphi u^{\f}): \Inc_{\mathcal{A}}(A^{\f}) \longrightarrow (A, f)$ exists a unique 1-cell in $\mathcal{A}$ of the form $p^{\varphi}:A^{\f}\longrightarrow \Alg_{\mathcal{A}}(B, h) = B^{\h}$  such that the following diagram commutes

\begin{equation*}\Label{1405161557}
\xy<1cm,0cm>
\POS (0,15) *+{\Inc_{\mathcal{A}}(A^{\f})} = "a11",
\POS (50,15) *+{\Inc_{\mathcal{A}} (B^{\h})} = "a13"
\POS (25,0) *+{(A, f)} = "a22"

\POS "a11" \ar^{\Inc_{\mathcal{A}}(\separ p^{\varphi})} "a13",
\POS "a11" \ar_{(pu^{\f}\!, \separ p\chi^{\f}\circ\separ \varphi u^{\f})\phantom{12}} "a22",
\POS "a13" \ar^{(u^{\h},\chi^{\h})} "a22",
\endxy
\end{equation*}\\

In particular, $pu^{\f} = u^{\h}p^{\varphi}$ and $p\chi^{\f}\circ\separ \varphi u^{\f} = \chi^{\h}p^{\varphi}$. Observe that the associated mate, to the first equality, is $\lambda = \varepsilon^{\h}p^{\varphi}d^{\f}\circ d^{\h}p\separ\eta^{\f}$ and that  $u^{\h}\lambda = \varphi$.\\

Consider a 2-cell of monads, $\theta:(p, \varphi)\longrightarrow (q,\psi):(A, f)\longrightarrow (B, h)$. Because of the construction of algebras for $\mathcal{A}$, the 2-adjunction provides an isomorphism of categories, for $A$ in $\mathcal{A}$ and $(X, f)$ in $\Mnd(\mathcal{A})$,

\begin{equation*}\Label{1405251634}
Hom_{\mathcal{A}}(A, \Alg_{\mathcal{A}}(X,f))\cong Hom_{\Mnd(\mathcal{A})}(\Inc_{\mathcal{A}}(A), (X,f))
\end{equation*}\\

\noindent given by the following assignment 

\begin{equation}\Label{1405251638}
\begin{array}{ccc}
\xy<1cm,0cm>
\POS (0,0) *+{A} = "a11",
\POS (25,0) *+{X^{\f}} = "a12",
\POS (9, 3.5) = "a",
\POS (9,-3.5) = "b"

\POS "a11" \ar@/^1.4pc/^{a} "a12",
\POS "a11" \ar@/_1.4pc/_{b} "a12",

\POS "a" \ar^{\phantom{1}\alpha} "b",

\endxy & \longmapsto &
\xy<1cm,0cm>
\POS (0,0) *+{(A, 1_{A})} = "a11",
\POS (30,0) *+{(X, f)} = "a13",
\POS (9, 3.5) = "a",
\POS (9,-3.5) = "b",

\POS (15, 0.5) *+{\scriptstyle{u^{\f}\alpha}} = "a12",

\POS "a11" \ar@/^1.5pc/^{(u^{\f}a,\separ \chi^{\f}a)} "a13",
\POS "a11" \ar@/_1.5pc/_{(u^{\f}b, \separ \chi^{\f}b)} "a13",

\POS "a" \ar "b",

\endxy
\end{array}
\end{equation}\\

\noindent cf. \cite{stro_fotm}. On the other hand, we have an equality of 2-cells\\

\begin{equation*}\Label{1405251719}
\begin{array}{ccc}
\xy<1cm,0cm>
\POS (0,0) *+{(A^{\f}, 1_{A^{\f}})} = "a11",
\POS (30,0) *+{(B, h)} = "a13",
\POS (10, 3.5) = "a",
\POS (10,-3.5) = "b",

\POS (15, 0.5) *+{\scriptstyle{\theta u^{\f}}} = "a12",

\POS "a11" \ar@/^1.5pc/^{(pu^{\f}, \separ p\chi^{\f}\circ\separ\varphi u^{\f})} "a13",
\POS "a11" \ar@/_1.5pc/_{(qu^{\f}, \separ q\chi^{\f}\circ\separ\psi u^{\f})} "a13",

\POS "a" \ar "b",

\endxy & = &
\xy<1cm,0cm>

\POS (0,0) *+{(A^{\f}, 1_{A^{\f}})} = "a11",
\POS (30,0) *+{(B, h)} = "a13",
\POS (11, 3.5) = "a",
\POS (11,-3.5) = "b",

\POS (16, 0.5) *+{\scriptstyle{\theta u^{\f}}} = "a12",

\POS "a11" \ar@/^1.5pc/^{(u^{\h}p^{\varphi}, \separ \chi^{\h}p^{\varphi})} "a13",
\POS "a11" \ar@/_1.5pc/_{(u^{\h}q^{\psi}, \separ \chi^{\h}q^{\psi})} "a13",

\POS "a" \ar "b",

\endxy
\end{array}
\end{equation*}\\

Therefore, to the 2-cell $\theta u^{\f}$ there corresponds, through the
assignment \eqref{1405251638}, a 2-cell $\Alg_{\mathcal{A}}(\theta
u^{\f})\eta^{AI}(A^{\f}):= \beta^{\theta}$, where
$\beta^{\theta}:p^{\varphi}\longrightarrow q^{\psi}$ and such that
$u^{\h}\beta^{\theta} = \theta u^{\f}$. We change the notation as follows
$\beta^{\theta} \separ\separ\separ\separ = \separ\separ\separ \separ\separ\widehat{\theta}$.\\

With these calculations at hand, we define the 2-functor $\Psi_{\E}$.

\begin{enumerate}

\item [1.-] On 0-cells, $(A, f)$ , $\Psi_{\E}(A, f) = d^{\f}\dashv u^{\f}$,
  \emph{i.e.} the Eilenberg-Moore adjunction.

\item [2.-] On 1-cells, $(p, \varphi):(A, f)\longrightarrow (B, h)$,  $\Psi_{\E}(p, \varphi) = (p, p^{\varphi}): d^{\f}\dashv u^{\f}\longrightarrow  d^{\h}\dashv u^{\h}$.

\item [3.-] On 2-cells, $\theta:(p, \varphi)\longrightarrow (q,\psi):(A,
  f)\longrightarrow (B, h)$, $\Psi_{\E}(\theta) = (\theta, \widehat{\theta}):(p, p^{\varphi})\longrightarrow (q, q^{\psi}): d^{\f}\dashv u^{\f}\longrightarrow  d^{\h}\dashv u^{\h}$.

\end{enumerate}

The unit and the counit for this 2-adjunction are given as follows. The component of the unit, at $l\dashv r$, is 
$\eta^{\Psi\Phi_{\!\E}}(l\dashv r):l\dashv
r\longrightarrow\Psi_{\E}\Phi_{\E}(l\dashv r)$, where
$\Psi_{\E}\Phi_{\E}(l\dashv r) = d^{\separ rl}\dashv u^{rl}$. \\

In \cite{stro_fotm}, Theorem 3, the author proved the existence of a \emph{comparison} 1-cell $k^{rl}:B\longrightarrow A^{rl}$, such that $u^{rl}k^{rl} = r$ and $d^{rl} = k^{rl}l$. Therefore, we can make the following definition $\eta^{\Psi\Phi_{\!\E}}(l\dashv r) = (1_{A}, k^{rl}, 1_{d^{rl}}): l\dashv r \longrightarrow d^{\separ rl}\dashv u^{rl}$.\\

In turn, the component of the counit, at $(A,f)$, is $\varepsilon^{\Phi\Psi_{\E}}(A, f):\Phi_{\E}\Psi_{\E}(A, f)\longrightarrow (A,f)$, where $\Phi_{\E}\Psi_{\E}(A, f) = (A, f)$. In this case, the counit is defined as $\varepsilon^{\Phi\Psi_{\E}}(A, f) = (1_{A}, 1_{\f}):(A, f)\longrightarrow (A, f)$. 

\newtheorem{1405251557}{Theorem}[section]
\begin{1405251557}
There exists a 2-adjunction $\Phi_{\E}\dashv \Psi_{\E}$. 
\end{1405251557}

\noindent \emph{Proof}:\\

We prove only one of the triangular identities and the other one is left to the reader. Using the definition of the unit and counit for this 2-adjunction, the triangular identity $\varepsilon^{\Phi\Psi_{\!\E}}\Phi_{\!\E}\circ\Phi_{\!\E}\eta^{\Psi\Phi_{\!\E}} = 1_{\Phi_{\E}}$ is proved as indicated.

\begin{eqnarray*}
(\varepsilon^{\Phi\Psi_{\!\E}}\Phi_{\!\E}\circ\Phi_{\!\E}\eta^{\Psi\Phi_{\!\E}})
(l\dashv r) &=& \varepsilon^{\Phi\Psi_{\!\E}}\Phi_{\!\E}(l\dashv r)\cdot\Phi_{\!\E}\eta^{\Psi\Phi_{\!\E}}
(l\dashv r) \\
&=& \varepsilon^{\Phi\Psi_{\!\E}}(A,rl)\cdot \Phi_{\!\E}(1_{A}, k^{rl},
1_{d^{rl}}) = (1_{A}, 1_{rl})\cdot(1_{A}, u^{rl}1_{d^{rl}})\\ 
&=& (1_{A}, 1_{rl}) = 1_{(A, rl)} = 1_{\Phi_{\E}(l\dashv r)} = 1_{\Phi_{\E}}(l\dashv r). 
\end{eqnarray*}

\begin{flushright}
$\square$
\end{flushright}

Since the right 2-adjoint assigns the Eilenberg-Moore adjunction to a monad
$(A,f)$, this 2-adjunction is called the \emph{Eilenberg-Moore 2-adjunction}.

\section{Eilenberg-Moore 2-Adjunction}

In this section, we apply the results of the Section 3 to the 2-category
$_{2}Cat$, the 2-category of small categories and functors, due to the fact
that this 2-category admits the construction of algebras. Therefore, we have a 2-adjunction

\begin{equation*}\Label{1310282200}
\begin{array}{cc}
\xymatrix@C=1.6cm{
\Adj_{\R}(_{2}Cat)\ar@<-3pt>[r]_-{\Phi_{\E}} & \Mnd\ar@<-3pt>[l]_-{\Psi_{\!\E}}(_{2}Cat) 
} &\!\!\! ,
\end{array}
\end{equation*}

Since the complete description, for a general $\mathcal{A}$, has been given above, we only give some remarks on the derived properties for this particular 2-category. \\

The description of the $2$-functor $\Psi_{\E}$, for this particular 2-category, is given by the following entries

\begin{enumerate}

\item [1.-] On $0$-cells, $\Psi_{\E}(\mathcal{C}, F) = D^{\F}\dashv U^{\F}$, \emph{i.e.} the Eilenberg-Moore adjunction.

\item [2.-] On $1$-cells, $(P,\varphi) : (\mathcal{C}, F)\longrightarrow (\mathcal{D}, H)$, $\Psi_{\E}(P, \varphi) = (P, P^{\separ\varphi}, \lambda^{\varphi})$. The action of the functor $P^{\separ\varphi}:\mathcal{C}^{\F}\longrightarrow \mathcal{D}^{\Hh}$ is the following

\begin{enumerate}

\item [(i)] On objects, $(M, \chi_{\M})$ in $\mathcal{C}^{\F}$, $P^{\separ\varphi}(M, \chi_{\M}) = (PM, P\chi_{\M}\cdot \varphi_{M})$.

\item [(ii)] On morphisms, $p$, $P^{\separ\varphi}(p) = Pp$.

\item [(iii)] The natural transformation $\lambda^{\varphi}$ is the mate of the identity $U^{\Hh}P^{\varphi} = P U^{\F}$. Using \eqref{1405151810}, we get the component of $\lambda^{\varphi}$ at $A$, in $\mathcal{C}$,

\begin{eqnarray*}
\lambda^{\varphi}A &=& \big(\varepsilon^{\D\!\U^{\!\Hh}}P^{\varphi}D^{\F}\circ D^{\Hh}\!P\eta^{\U\!\D^{\!\F}}\big)(A)\\
&=& P\mu^{\F}A\cdot \varphi FA\cdot HP\eta^{\F}A\\
&=& \varphi A.
\end{eqnarray*}

\end{enumerate}

\item [3.-] On $2$-cells, $\theta:(P,\varphi)\longrightarrow(Q,\psi)$, we have

\begin{equation*}\Label{1403152054}
\Psi_{\E}(\theta) = (\alpha^{\theta}, \beta^{\theta}) = (\theta, \hat{\theta})
\end{equation*}\\

The induced natural tranformation $\hat{\theta}:P^{\varphi}\longrightarrow Q^{\psi}:\mathcal{C}^{\F}\longrightarrow \mathcal{D}^{\Hh}$ is defined through its components, using the condition $ \theta\separ U^{\F}=U^{\Hh}\hat{\theta}$, as

\begin{equation*}\Label{1403152055}
\hat{\theta}(M,\chi_{\M}) = \theta M.
\end{equation*}

\end{enumerate}

Since we have a 2-adjunction, the following isomorphism of categories takes place, natural for $L\dashv R$ in $\Adj_{\R}({}_{2}Cat)$ and $(\mathcal{X}, H)$ in $\Mnd({}_{2}Cat)$:

\begin{equation}\Label{1403191902}
Hom_{\Adj_{\R}({}_{2}Cat)}(L\dashv R, \Psi_{\E}(\mathcal{X}, H)) \cong Hom_{\Mnd({}_{2}Cat)}(\Phi_{\E}(L\dashv R), (\mathcal{X}, H))
\end{equation}


\section{Monoidal Liftings (Eilenberg-Moore Type)} 

\subsection{Colax Monads}

In this section, we give the definition of a colax monad.\\

\newtheorem{1403161259}{Definition}[section]
\begin{1403161259}
A colax monad $\big((F, \xi, \gamma), \mu^{\F}, \eta^{\F}\big)$ over the monoidal category $(\mathcal{C}, \otimes, I)$ consists of the following

\begin{enumerate}

\item $(F,\mu^{\F}, \eta^{\F})$ is a monad on $\mathcal{C}$.

\item $(F,\xi,\gamma):(\mathcal{C}, \otimes, I)\longrightarrow (\mathcal{C}, \otimes, I)$ is a colax monoidal functor. That is to say, the natural transformations $\xi:F\cdot\otimes\longrightarrow\otimes\cdot (F\times F) $ and $\gamma:F\cdot\delta_{I}\longrightarrow \delta_{I}$  fulfills the commutativity of the following diagrams

\begin{equation}\Label{1403181443}
\xymatrix@C=1.8cm@R=1.4cm{
F((A\otimes B)\otimes C)\ar[r]^{\xi_{A\otimes B,C}}\ar[d]_{Fa_{A,B,C}} & F(A\otimes B)\otimes FC\ar[r]^{\xi_{A,B}\otimes\: FC} & (FA\otimes FB)\otimes FC\ar[d]^{a_{FA,FB,FC}}\\
F(A\otimes (B\otimes C))\ar[r]_{\xi_{A, B\otimes C}} & FA\otimes F(B\otimes C)\ar[r]_{FA\:\otimes\:\xi_{B,C}} & FA\otimes (FB\otimes FC)
}
\end{equation}

\begin{equation}\Label{1403191813}
\begin{array}{ccc}
\xymatrix@C=1.2cm@R=1.2cm{
F(I\otimes A)\ar[r]^-{\xi_{I, A}}\ar@/_0.8pc/[rd]_{Fl_{A}} & FI \otimes FA\ar[r]^-{\gamma\separ\otimes FA} & I\otimes FA\ar@/^0.8pc/[ld]^{l_{FA}} \\
 & FA  & 
} &
&
\xymatrix@C=1.2cm@R=1.2cm{
FA\otimes I\ar@/_0.8pc/[rd]_{r_{FA}} & FA\otimes FI\ar[l]_-{FA\separ \otimes \gamma} & F(A\otimes I)\ar[l]_-{\xi_{A,I}}\ar@/^0.8pc/[ld]^{Fr_{A}}\\
 & FA & 
}
\end{array}
\end{equation}

\item $\mu^{\F}:(F,\xi,\gamma)\cdot(F,\xi,\gamma)\longrightarrow
  (F,\xi,\gamma)$ and $\eta^{\F}:(1_{\mathcal{C}}, 1_{\otimes},
  1_{\delta_{I}})\longrightarrow (F,\xi,\gamma)$ are colax natural
  transformations, \emph{i.e.} apart from the fact that they are natural
  transformations, they have to fulfill the following commutative diagrams

\begin{equation}\Label{1403191825}
\begin{array}{ccc}
\xymatrix@C=1.5cm@R=1.5cm{
FF\otimes \ar[r]^-{F\xi}\ar[d]_{\mu^{\F}\otimes} & F\otimes (F\times F)\ar[r]^-{\xi (F\times F) } & \otimes (FF\times FF)\ar[d]^{\otimes(\mu^{\F}\times \mu^{\F})}\\
F\otimes\ar[rr]_{\xi} & &  \otimes (F\times F)   
} & 
\xy<1cm,0cm>
\POS (0,-20.5) *+{,},
\endxy
&
\xymatrix@C=1.15cm@R=1.5cm{
FF\delta_{I}\ar[d]_{\mu^{\F}\delta_{I}}\ar[r]^-{F\gamma}  & F\delta_{I}\ar[r]^-{\gamma} & \delta_{I}\\
F\delta_{I}\ar@/_1pc/[rru]_{\gamma} & &     
}
\end{array}
\end{equation}

\noindent \&

\begin{equation}\Label{1403191842}
\begin{array}{ccc}
\xy
\POS (-22, 0) *+{\otimes} = "a",
\POS (22, 0) *+{\otimes} = "d",
\POS (-22, -20) *+{F\otimes} = "c", 
\POS (22, -20) *+{\otimes(F\times F)} = "b",
\POS "a" \ar^{1_{\otimes}} "d",
\POS "d" \ar^{\otimes(\eta^{\F}\times\eta^{\F})} "b",
\POS "a" \ar_{\eta^{\F}\otimes} "c",
\POS "c" \ar_{\xi} "b",
\endxy
& 
\xy<1cm,0cm>
\POS (0,-20.5) *+{,},
\endxy
&
\xy<1cm,0cm>
\POS (0, -4)
\xymatrix@C=1.2cm@R=1cm{
\delta_{I}\ar[d]_{\eta^{\F}\delta_{I}} \ar[r]^{1_{\delta_{I}}} &  \delta_{I}\\
F\delta_{I}\ar@/_0.7pc/[ru]_{\gamma} &     
}
\endxy
\end{array}
\end{equation}

\end{enumerate}

\end{1403161259}

Since the natural transformation $\gamma$ has only one component, at $0$ in $\mathbf{1}$, then this natural transformation and its component will be denoted indistinctly as $\gamma$.\\


Using the isomorphism \eqref{1403191902}, the following bijection can be obtained, cf. \cite{moie_motc}

\newtheorem{1403192127}{Theorem}[section]
\begin{1403192127}

There is bijective correspondance between the following structures

\begin{enumerate}

\item [1.-] Colax monads $\big((F, \xi, \gamma), \mu^{\F}, \eta^{\F}\big)$, for the monoidal structure $(\mathcal{C}, \otimes, I, a, l, r)$.

\item [2.-] Morphisms and natural transformations of monads of the form 

\begin{eqnarray*}
(\otimes, \xi)&:& (\mathcal{C}\times\mathcal{C}, F\times F)\longrightarrow (\mathcal{C}, F),\\
(\delta_{I}, \gamma)&:&(\mathbf{1}\:, 1_{\mathbf{1}}\:)\longrightarrow (\mathcal{C},F)\\
a&:&(\otimes\cdot(\otimes\times \mathcal{C}),\otimes(\xi\times F)\circ\xi(\otimes\times\mathcal{C}))\longrightarrow(\otimes\cdot (\mathcal{C}\times\otimes)\cdot a_{\mathcal{C}}, \otimes(F\times\xi)a_{\mathcal{C}}\circ\xi(\mathcal{C}\times\otimes)a_{\mathcal{C}})\\
&&\quad\quad\quad\quad\quad\quad\quad\quad\quad\quad\quad\quad:((\mathcal{C}\times\mathcal{C})\times\mathcal{C}, (F\times F)\times F)\longrightarrow (\mathcal{C}, F), \\
l&:& (\otimes\cdot(\delta_{I}\times \mathcal{C})\cdot l_{\mathcal{C}}^{-1}, \otimes(\gamma\times F)\separ l^{-1}_{\mathcal{C}}\circ
  \separ\xi(\delta_{I}\times\mathcal{C})\separ l^{-1}_{\mathcal{C}})\longrightarrow (1_{\mathcal{C}}, 1_{F}):(\mathcal{C}, F)\longrightarrow (\mathcal{C}, F),\\
r&:&(\otimes\cdot( \mathcal{C}\times\delta_{I})\cdot r_{\mathcal{C}}^{-1},
\otimes(F\times\gamma)\separ r^{-1}_{\mathcal{C}}\circ
\xi(\mathcal{C}\times\delta_{I})\separ r^{-1}_{\mathcal{C}})\longrightarrow (1_{\mathcal{C}}, 1_{F}):(\mathcal{C}, F)\longrightarrow (\mathcal{C}, F).\\ 
\end{eqnarray*}

\item [3.-] Monoidal structures for the Eilenberg-Moore category, $(\mathcal{C}^{\F}, \widehat{\otimes} , \hat{I}, \hat{a}, \hat{l}, \hat{r})$  such that the following diagram of arrows and surfaces commutes

\begin{equation}\Label{1303192133}
\begin{array}{cc}
\xy<1cm,0cm>
\POS (0, 10) *+{(a)},
\POS (0, 0) *+{\mathcal{C}\times\mathcal{C}} = "c1",
\POS (25, 0) *+{\mathcal{C}} = "c2", 
\POS (0, -20) *+{\mathcal{C}^{\F}\times\mathcal{C}^{\F}} = "d1",
\POS (25, -20) *+{\mathcal{C}^{\F}} = "d2", 
\POS (18, -10) = "a",
\POS (10, -18) = "b",
\POS "c1" \ar^-{\otimes} "c2",
\POS "d1" \ar^{U^{\F}\times U^{\F}} "c1",
\POS "d2" \ar_{U^{\F}} "c2",
\POS "d1" \ar_-{\widehat{\otimes}} "d2",
\endxy & 
\xy<1cm,0cm>
\POS (0, 10) *+{(b)},
\POS (0,0) *+{\mathbf{1}} = "c1",
\POS (25,0) *+{\mathcal{C}} = "c2",
\POS (25,-20) *+{\mathcal{C}^{\F}} = "d2",
\POS (0,-20) *+{\mathbf{1}^{1_{\mathbf{1}}}} = "d1",
\POS "c1" \ar^-{\delta_{I}} "c2",
\POS "d1" \ar^{U^{1_{\mathbf{1}}}} "c1",
\POS "d2" \ar_{U^{\F}} "c2",
\POS "d1" \ar_-{\delta_{\hat{I}}} "d2",
\endxy
\end{array}
\end{equation}

\begin{equation*}\Label{1404271827}
\begin{array}{ccc}
\xy<1cm,0cm>
\POS (0, 0) *+{\mathcal{C}^{3}} = "c1",
\POS (7, 0) *+{\phantom{(\mathcal{C})}}= "c3",
\POS (30, 0) *+{\mathcal{C}} = "c2", 
\POS (11, 4) = "a",
\POS (11, -4) = "b",
\POS (11,-36) = "c",
\POS (11,-44) = "d",

\POS (0,-40) *+{(\mathcal{C}^{\F})^{3}} = "d1",
\POS (30,-40) *+{\mathcal{C}^{\F}} = "d2",

\POS "c1" \ar@/^1.5pc/^{\otimes\cdot (\otimes\times \mathcal{C})} "c2",
\POS "c1" \ar@/_1.5pc/_{\otimes\cdot (\mathcal{C}\times\otimes)\cdot a_{\mathcal{C}}} "c2",
\POS "a" \ar^{\phantom{11}a} "b",

\POS "d1" \ar@/^1.5pc/^{\widehat{\otimes}\cdot (\widehat{\otimes}\times \mathcal{C}^{\F})} "d2",
\POS "d1" \ar@/_1.5pc/_{\widehat{\otimes}\cdot (\mathcal{C}^{\F}\times\widehat{\otimes})\cdot a_{\mathcal{C}^{\F}}} "d2",
\POS "c" \ar^{\phantom{11}\widehat{a}} "d",

\POS "d1" \ar^{(U^{\F})^{3}} "c1",
\POS "d2" \ar_{U^{\F}} "c2",
\endxy &
\xy<1cm,0cm>
\POS (0, 0) *+{\mathcal{C}} = "c1",
\POS (7, 0) *+{\phantom{(\mathcal{C})}}= "c3",
\POS (30, 0) *+{\mathcal{C}} = "c2", 
\POS (11, 4) = "a",
\POS (11, -4) = "b",
\POS (11,-36) = "c",
\POS (11,-44) = "d",

\POS (0,-40) *+{\mathcal{C}^{\F}} = "d1",
\POS (30,-40) *+{\mathcal{C}^{\F}} = "d2",

\POS "c1" \ar@/^1.5pc/^{\otimes\cdot (\delta_{I}\times \mathcal{C})\cdot l^{-1}_{\mathcal{C}}} "c2",
\POS "c1" \ar@/_1.5pc/_{1_{\mathcal{C}}} "c2",
\POS "a" \ar^{\phantom{12}l} "b",

\POS "d1" \ar@/^1.5pc/^{\widehat{\otimes}\cdot (\delta_{\hat{I}}\times
  \mathcal{C}^{\F})\cdot l^{-1}_{\mathcal{C}^{\F}}} "d2",
\POS "d1" \ar@/_1.5pc/_{1_{\mathcal{C}^{\F}}} "d2",
\POS "c" \ar^{\phantom{12}\widehat{l}} "d",

\POS "d1" \ar^{U^{\F}} "c1",
\POS "d2" \ar_{U^{\F}} "c2",
\endxy
&
\xy<1cm,0cm>
\POS (0, 0) *+{\mathcal{C}} = "c1",
\POS (7, 0) *+{\phantom{(\mathcal{C})}}= "c3",
\POS (30, 0) *+{\mathcal{C}} = "c2", 
\POS (11, 4) = "a",
\POS (11, -4) = "b",
\POS (11,-36) = "c",
\POS (11,-44) = "d",

\POS (0,-40) *+{\mathcal{C}^{\F}} = "d1",
\POS (30,-40) *+{\mathcal{C}^{\F}} = "d2",

\POS "c1" \ar@/^1.5pc/^{\otimes\cdot (\mathcal{C}\times\delta_{I})\cdot r^{-1}_{\mathcal{C}}} "c2",
\POS "c1" \ar@/_1.5pc/_{1_{\mathcal{C}}} "c2",
\POS "a" \ar^{\phantom{12}r} "b",

\POS "d1" \ar@/^1.5pc/^{\widehat{\otimes}\cdot
  (\mathcal{C}^{\F}\times\delta_{\hat{I}})\cdot r^{-1}_{\mathcal{C}^{\F}}} "d2",
\POS "d1" \ar@/_1.5pc/_{1_{\mathcal{C}^{\F}}} "d2",
\POS "c" \ar^{\phantom{12}\widehat{r}} "d",

\POS "d1" \ar^{U^{\F}} "c1",
\POS "d2" \ar_{U^{\F}} "c2",
\endxy
\end{array}
\end{equation*}

\end{enumerate}

\end{1403192127}

\noindent\emph{Proof}:\\

$1\Rightarrow 2$)\\

Consider a colax monad $\big((F, \xi, \gamma), \mu^{\F}, \eta^{\F}\big)$, for the monoidal structure $(\mathcal{C}, \otimes, I)$. In particular, the multiplication and the unit of the monad are colax natural transformations and the first diagrams in \eqref{1403191825} and \eqref{1403191842} commute. Therefore, we have a monad morphism $(\otimes, \xi):(\mathcal{C}\times\mathcal{C} , F\times F)\longrightarrow (\mathcal{C}, F)$.\\

Likewise, the commutativity of the second diagrams in \eqref{1403191825} and \eqref{1403191842} implies that $(\delta_{I}, \gamma):(\mathbf{1}\:, 1_{\mathbf{1}}\:)\longrightarrow (\mathcal{C},F)$ is a morphism of monads. Note that the requirement $(\delta_{I}, \gamma)$ is a monad morphism is equivalent to the statement $(I, \gamma)$ is an Eilenberg-Moore algebra.\\

Since $(\otimes, \xi)$ is a morphism of monads then the following are also morphisms of monads $\big(\otimes\cdot(\otimes\times\mathcal{C}), \otimes(\xi\times F)\circ\xi(\otimes\times\mathcal{C})\big)$ and $(\otimes\cdot(\mathcal{C}\times\otimes)\cdot a_{\mathcal{C}}, \otimes(F\times\xi)a_{\mathcal{C}}\circ \xi(\mathcal{C}\times\otimes)a_{\mathcal{C}})$ from $((\mathcal{C}\times \mathcal{C})\times \mathcal{C}, (F\times F)\times F)$ to $(\mathcal{C}, F)$. Furthermore, due to the commutativity of the diagram \eqref{1403181443}, the following is a $2$-cell in $\Mnd({}_{2}Cat)$

\begin{equation*}
\xy<1cm,0cm>
\POS (0, 0) *+{((\mathcal{C}\times\mathcal{C})\times\mathcal{C}, (F\times F)\times F)} = "c1",
\POS (7, 0) *+{\phantom{(\mathcal{C})}}= "c3",
\POS (50, 0) *+{(\mathcal{C}, F)} = "c2", 
\POS (27, 5) = "a",
\POS (27, -5) = "b",
\POS "c3" \ar@/^2pc/^{(\otimes\cdot (\otimes\times \mathcal{C}),\:\otimes(\xi\times F)\circ\xi(\otimes\times\mathcal{C}))} "c2",
\POS "c3" \ar@/_2pc/_{(\otimes\cdot (\mathcal{C}\times\otimes)\cdot a_{\mathcal{C}},\:\otimes(F\times\xi)a_{\mathcal{C}}\circ\xi(\mathcal{C}\times\otimes)a_{\mathcal{C}})} "c2",
\POS "a" \ar^{\phantom{1}a} "b",
\endxy
\end{equation*}\\

Likewise, because $(\otimes, \xi)$ and $(\delta_{I}, \gamma)$ are monad morphisms, $(\otimes\cdot (\delta_{I}\times \mathcal{C})\cdot l^{-1}_{\mathcal{C}}, \otimes(\gamma\times F)\separ l^{-1}_{\mathcal{C}}\circ   \separ\xi(\delta_{I}\times\mathcal{C})\separ l^{-1}_{\mathcal{C}})$ is also  a monad morphism. Using the commutativity of the first diagram in \eqref{1403191813}, we can consider the monad $2$-cell

\begin{equation*}
\xy<1cm,0cm>
\POS (0, 0) *+{(\mathcal{C}, F)} = "c1",
\POS (7, 0) *+{\phantom{(\mathcal{C})}}= "c3",
\POS (55, 0) *+{(\mathcal{C}, F)} = "c2", 
\POS (19, 5) = "a",
\POS (19, -5) = "b",
\POS "c1" \ar@/^2pc/^{(\otimes\cdot (\delta_{I}\times \mathcal{C})\separ
  l^{-1}_{\mathcal{C}},\:\otimes\!\phantom{.}(\gamma\times
  F)\!\phantom{.}l^{-1}_{\mathcal{C}}\circ\separ \xi (\delta_{I}\times \mathcal{C})l^{-1}_{\mathcal{C}})} "c2",
\POS "c1" \ar@/_2pc/_{(1_{\mathcal{C}}, 1_{\F})} "c2",
\POS "a" \ar^{\phantom{1234}l} "b",
\endxy
\end{equation*}\\

In a similar way, the following is a monad transformation, $r:(\otimes\cdot(\mathcal{C}\times\delta_{I})\cdot r_{\mathcal{C}}^{-1}, \otimes(F\times\gamma)\separ r^{-1}_{\mathcal{C}}\circ\xi(\mathcal{C}\times\delta_{I})\separ r^{-1}_{\mathcal{C}})\longrightarrow (1_{\mathcal{C}}, 1_{F}):(\mathcal{C}, F)\longrightarrow (\mathcal{C}, F)$.\\

$2\Rightarrow 1$)\\ 

Note that the aforementioned claims can be reverted.\\

2 $\Rightarrow$ 3)\\

Take the monad morphism $(\otimes, \xi): (\mathcal{C}\times\mathcal{C},F\times F)\longrightarrow (\mathcal{C}, F)$. In order to use the isomorphism \eqref{1403191902}, we make $L\dashv R = D^{\F}\times D^{\F} \dashv U^{\F}\times U^{\F}$ and $(\mathcal{X}, H,\mu^{\Hh},\eta^{\Hh}) = (\mathcal{C}, F, \mu^{\F}, \eta^{\F})$. Therefore, to this monad morphism corresponds a morphism of adjunctions of the form $(\otimes, \otimes^{\xi}): D^{\F}\times D^{\F} \dashv U^{\F}\times U^{\F} \longrightarrow D^{\F}\dashv U^{\F}$ such that a diagram like (\ref{1303192133}a) commutes. According to the definition of $\Psi^{\E}$, the functor $\otimes^{\xi}$ acts, on objects, as follows

\begin{equation*}\Label{1403201159}
\otimes^{\xi}\big((M,\chi_{\M}),(N,\chi_{\N})\big) = \big(\otimes\! (M, N), \otimes(\chi_{\M}, \chi_{\N})\cdot \xi_{\M,\N}\big)
\end{equation*}\\

The previous action is defined at the beginning of the proof of Theorem 7.1, \cite{moie_motc}. On morphisms, we have

\begin{equation*}\Label{1403201301}
\otimes^{\xi}(p,q) = \otimes(p,q)
\end{equation*}\\

\noindent We change the notation from $\otimes^{\xi}$ to $\widehat{\otimes}$.\\

If in the isomorphism \eqref{1403191902}, we  make $L\dashv R = 1_{\mathbf{1}}\dashv 1_{\mathbf{1}}$ and $(\mathcal{X}, H,\mu^{\Hh},\eta^{\Hh}) = (\mathcal{C}, F, \mu^{\F}, \eta^{\F})$. The monad morphism $(\delta_{I}, \gamma)$ has an associated morphism of adjunctions of the form $(\delta_{I}, \delta_{I}^{\phantom{I}\gamma}):(1_{\mathbf{1}}\dashv 1_{\mathbf{1}})\longrightarrow D^{\F}\dashv U^{\F}$ such that a diagram like (\ref{1303192133}b) commutes. According to the definition of $\Psi^{\E}$, the functor $\delta_{I}^{\phantom{I}\gamma}$ acts as follows

\begin{equation*}\Label{1403201409}
\delta_{I}^{\phantom{I}\gamma}(0, 1_{0}) = (\delta_{I}(0), \delta_{I}(1_{0})\cdot \gamma) = (I, \gamma). 
\end{equation*}\\

\noindent On morphisms,  

\begin{equation*}\Label{1403201416}
\delta_{I}^{\phantom{I}\gamma}(1_{0}) = \delta_{I}(1_{0}) = 1_{I} = 1_{(I,\gamma)}.
\end{equation*}\\

\noindent If we make the following definition $\hat{I} = (I, \gamma)$, then $\delta_{I}^{\phantom{I}\gamma} := \delta_{\hat{I}}$. The algebra $(I, \gamma)$ is the unit of the monoidal structure on $\mathcal{C}^{\F}$. \\

Suppose that we have a natural transformation of the form $a:(\otimes\cdot(\otimes\times \mathcal{C}),\otimes(\xi\times F)\circ\xi(\otimes\times\mathcal{C}))\longrightarrow(\otimes\cdot (\mathcal{C}\times\otimes)\cdot a_{\mathcal{C}}, \otimes(F\times\xi)a_{\mathcal{C}}\separ\circ\separ\xi(\mathcal{C}\times\otimes)a_{\mathcal{C}}):((\mathcal{C}\times\mathcal{C})\times\mathcal{C}, (F\times F)\times F)\longrightarrow (\mathcal{C}, F)$ then we can make $L\dashv R = (D^{\F}\times D^{\F})\times D^{\F}\dashv (U^{\F}\times U^{\F})\times U^{\F}$ and $(\mathcal{X}, H,\mu^{\Hh},\eta^{\Hh}) = (\mathcal{C}, F, \mu^{\F}, \eta^{\F})$. Therefore, to the previous $2$-cell of monads, we can associate a $2$-cell of adjunctions of the form 

\begin{equation*}
\xy<1cm,0cm>
\POS (0, 0) *+{\mathcal{C}^{3}} = "c1",
\POS (7, 0) *+{\phantom{(\mathcal{C})}}= "c3",
\POS (40, 0) *+{\mathcal{C}} = "c2", 
\POS (14, 4) = "a",
\POS (14, -4) = "b",
\POS (14,-36) = "c",
\POS (14,-44) = "d",

\POS (0,-40) *+{(\mathcal{C}^{\F})^{3}} = "d1",
\POS (40,-40) *+{\mathcal{C}^{\F}} = "d2",

\POS "c1" \ar@/^1.7pc/^{\otimes\cdot (\otimes\times \mathcal{C})} "c2",
\POS "c1" \ar@/_1.7pc/_{\otimes\cdot (\mathcal{C}\times\otimes)\cdot a_{\mathcal{C}}} "c2",
\POS "a" \ar^{\phantom{12}a} "b",

\POS "d1" \ar@/^1.7pc/^{[\otimes\cdot (\otimes\times \mathcal{C})]^{\separ{}^{\cdot}\!\xi^{2}}} "d2",
\POS "d1" \ar@/_1.7pc/_{[\otimes\cdot (\mathcal{C}\times\otimes)\cdot a_{\mathcal{C}}]^{{}_{\cdot}\xi^{2}}} "d2",
\POS "c" \ar^{\phantom{12}\beta^{a}} "d",

\POS "c1" \ar@<-4pt>_{(D^{\F})^{3}} "d1",
\POS "d1" \ar@<-1pt>_{(U^{\F})^{3}} "c1",
\POS "c2" \ar@<-2pt>_{D^{\F}} "d2",
\POS "d2" \ar@<-3pt>_{U^{\F}} "c2",
\endxy
\end{equation*}\\
 
In order to reduce expressions, we used and will be using the following notation 

\begin{eqnarray*}
{}^{\cdot}\xi^{2} &:=& \otimes(\xi\times F)\circ\xi(\otimes\times\mathcal{C}), \\
{}_{\cdot}\xi^{2} &:=& \otimes(F\times\xi)a_{\mathcal{C}}\circ\xi(\mathcal{C}\times\otimes)a_{\mathcal{C}}, \\
(\cdot)^{3} &:=& (\cdot\times\cdot)\times\cdot 
\end{eqnarray*}

Since $\Psi_{\E}$ is a 2-functor, we have

\begin{eqnarray*}
[\otimes\cdot (\otimes\times \mathcal{C})]^{\separ{}^{\cdot}\!\xi^{2}} &=& 
\Psi_{\E}(\otimes\cdot (\otimes\times \mathcal{C}), \otimes(\xi\times F)\circ\xi(\otimes\times\mathcal{C})
)\\
&=& \Psi_{\E}\big( (\otimes, \xi)\cdot (\otimes\times\mathcal{C}, \xi\times
F)\big)\\
&=& \Psi_{\E}(\otimes, \xi)\cdot \Psi_{\E}(\otimes\times\mathcal{C}, \xi\times
F)\\
&=& \otimes^{\xi}\cdot (\otimes\times\mathcal{C})^{\xi\times F} =
\widehat{\otimes}\cdot (\widehat{\otimes}\times\mathcal{C}^{\F})
\end{eqnarray*}

$\phantom{br}$

In the same way, we can check that
$[\otimes\cdot(\mathcal{C}\times\otimes)\cdot
a_{\mathcal{C}}]^{{}_{\cdot}\xi^{2}} = \widehat{\otimes}\cdot
(\mathcal{C}^{\F}\times\widehat{\otimes})\cdot a_{\mathcal{C}^{\F}}$. \\

We change the notation $\beta^{a}$ for $\widehat{a}$ and we get a natural
transformation 

\begin{equation*}
\widehat{a}:\widehat{\otimes}(\widehat{\otimes}\times\mathcal{C}^{\F})\longrightarrow
\widehat{\otimes}(\mathcal{C}^{\F}\times\widehat{\otimes})\cdot
a_{\mathcal{C}^{\F}}:(\mathcal{C}^{\F}\times\mathcal{C}^{\F})\times\mathcal{C}^{\F}\longrightarrow
\mathcal{C}^{\F}.
\end{equation*}

$\phantom{br}$ 

Using the definition of the functor $\Psi_{\E}$ on the $2$-cell $a$, we get
the components as 

\begin{equation*}\Label{1403221055}
\widehat{a}\big(\big((M,\chi_{\M}),(N,\chi_{\N})\big), (M',\chi_{\M'})\big)= a(M,N,M')
\end{equation*}\\

Suppose we have a $2$-cell in $\Mnd({}_{2}Cat)$ of the form $l:(\otimes\cdot(\delta_{I}\times \mathcal{C})\cdot l_{\mathcal{C}}^{-1}, \otimes(\gamma\times F)\separ l^{-1}_{\mathcal{C}}\circ \separ\xi(\delta_{I}\times\mathcal{C})\separ l^{-1}_{\mathcal{C}})\longrightarrow (1_{\mathcal{C}}, 1_{F}):(\mathcal{C}, F)\longrightarrow (\mathcal{C}, F)$. If in the isomorphism \eqref{1403191902}, we make $L\dashv R = D^{\F} \dashv U^{\F}$ and $(\mathcal{X},H,\mu^{\Hh},\eta^{\Hh}) = (\mathcal{C}, F, \mu^{\F}, \eta^{\F})$, it can be obtained a $2$-cell in the $2$-category $\Adj_{\R}({}_{2}Cat)$ of the form $(l,\beta^{l}):(\otimes\cdot (\delta_{I}\times \mathcal{C})\cdot l^{-1}_{\mathcal{C}}, 1_{\mathcal{C}})\longrightarrow (\lbrack \otimes\cdot (\delta_{I}\times \mathcal{C})\cdot l^{-1}_{\mathcal{C}}\rbrack^{\gamma\circ\xi}, \lbrack 1_{\mathcal{C}}\rbrack^{1_{\F}}):D^{\F}\dashv U^{\F}\longrightarrow D^{\F}\dashv U^{\F}$. Where we used the notation $\gamma\circ\xi = \otimes\!\phantom{.}(\gamma\times F)\!\phantom{.}l^{-1}_{\mathcal{C}}\circ\xi (\delta_{I}\times \mathcal{C})\separ l^{-1}_{\mathcal{C}}$. We change the notation from $\beta^{l}$ to $\hat{l}$.\\

In the same way as before, it can be proved that $\lbrack \otimes\cdot (\delta_{I}\times \mathcal{C})\cdot l^{-1}_{\mathcal{C}}\rbrack^{\gamma\circ\xi} = \widehat{\otimes}\separ (\delta_{\hat{I}}\times \mathcal{C}^{\F})\separ  l^{-1}_{\mathcal{C}^{\F}} $ and $\lbrack 1_{\mathcal{C}}\rbrack^{1_{\F}} =  1_{\mathcal{C}^{\F}}$. Therefore, we obtain a natural transformation   $\hat{l}:\widehat{\otimes}\separ(\delta_{\hat{I}}\times \mathcal{C}^{\F})\separ   l^{-1}_{\mathcal{C}^{\F}} \longrightarrow 1_{\mathcal{C}^{\F}}:\mathcal{C}^{\F}\longrightarrow \mathcal{C}^{\F}   $. Using the definition of the $2$-functor $\Psi_{E}$ on the $2$-cell $l$, the component of the natural transformation $\hat{l}$ on $(M, \chi_{\M})$ is

\begin{equation*}\Label{1405071231}
\hat{l}(M,\chi_{\M}) = l_{\M}
\end{equation*}\\

Similarly to the $2$-cell $r:(\otimes\cdot(\mathcal{C}\times\delta_{I})\cdot r^{-1}_{\mathcal{C}}, \:\otimes\!\phantom{.}(F\times \gamma)\!\phantom{.}r^{-1}_{\mathcal{C}}\circ\xi\separ (\mathcal{C}\phantom{.}\!\times \delta_{I})\phantom{.}\! r^{-1}_{\mathcal{C}})\longrightarrow(1_{\mathcal{C}}, 1_{\F}):(\mathcal{C}, F)\longrightarrow(\mathcal{C}, F)$ there corresponds a natural transformation $\hat{r}:\widehat{\otimes} \separ(\mathcal{C}^{\F}\times \delta_{\hat{I}})\separ r^{-1}_{\mathcal{C}^{\F}}\longrightarrow 1_{\mathcal{C}^{\F}}: \mathcal{C}^{\F}\longrightarrow \mathcal{C}^{\F}$. The component of this natural transformation, at $(M,\chi_{\M})$, is 

\begin{equation}\Label{1405071326}
\hat{r}(M,\chi_{\M}) = r_{\M}
\end{equation}\\

Since the natural transformations $a, l$ and $r$ fulfill the coherence conditions for a monoidal struture and $U^{\F}$ is faithfull then $\hat{a}, \hat{l}$ and $\hat{r}$ fulfill the pentagon and the triangle coherence conditions. Therefore, $(\mathcal{C}^{\F}, \widehat{\otimes}, \hat{I}, \hat{a},\hat{l}, \hat{r})$ is a monoidal structure over $\mathcal{C}^{\F}$.\\

3 $\Rightarrow$ 2)\\

Note that the aforementioned statements can be reverted. For example, take the
morphism of adjunctions
$(a,\widehat{a}):(\otimes\cdot(\otimes\times\mathcal{C}),
\widehat{\otimes}\cdot(\widehat{\otimes}\times\mathcal{C}^{\F}))\longrightarrow
(\otimes\cdot(\mathcal{C}\times\otimes)\cdot a_{\mathcal{C}},
\widehat{\otimes}\cdot(\mathcal{C}^{\F}\times\widehat{\otimes})\cdot
a_{\mathcal{C}^{\F}}):(U^{\F}\times U^{\F})\times U^{\F}\dashv (D^{\F}\times
D^{\F})\times D^{\F}\longrightarrow U^{\F}\dashv D^{\F}$. The image of this
2-cell, under $\Phi_{\E}$, is  $a:(\otimes, \varphi_{\otimes})\cdot(\otimes\times\mathcal{C}, \varphi_{\otimes}\times F)\longrightarrow(\otimes,\varphi_{\otimes})\cdot(\mathcal{C}\times\otimes, F\times\varphi_{\otimes})\cdot(a_{\mathcal{C}}, 1_{F\times(F\times F)\cdot a_{\mathcal{C}}}):((\mathcal{C}\times\mathcal{C})\times\mathcal{C}, (F\times F)\times F)\longrightarrow (\mathcal{C}, F)$, \emph{i.e.}

\begin{equation*}
a:(\otimes\cdot(\otimes\times\mathcal{C}), \otimes(\varphi_{\otimes}\times
F)\circ\varphi_{\otimes}(\otimes\times\mathcal{C}))\longrightarrow
(\otimes\cdot(\mathcal{C}\times\otimes)\cdot a_{\mathcal{C}},
\otimes(F\times\varphi_{\otimes})\cdot a_{\mathcal{C}}\circ
\varphi_{\otimes}\cdot(\mathcal{C}\times\otimes)\cdot
a_{\mathcal{C}}):(\mathcal{C}^{3}, F^{3})\longrightarrow (\mathcal{C}, F)
\end{equation*}

\noindent is a 2-cell in $\Mnd({}_{2} Cat)$.\\

Everytime we used the isomorphism \eqref{1403191902}, the monad $(\mathcal{C}, F, \mu^{\F}, \eta^{\F})$ was always
taken fixed, therefore the implication 2 $\Rightarrow$ 3 is natural in the monad $(\mathcal{C}, F, \mu^{\F}, \eta^{\F})$. 

\begin{flushright}
$\square$
\end{flushright}


\section{Kleisli 2-Adjunction}

Based on either \cite{brto_eicb} or \cite{cljs_klem}, the following 2-adjunction takes place

\begin{equation*}\Label{1310282201}
\begin{array}{cc}
\xymatrix@C=1.6cm{
\Mnd^{\scriptscriptstyle{\bullet}}({}_{2}Cat) \ar@<-3pt>[r]_-{\Psi_{\K}} & \Adj_{\Ll}({}_{2}Cat)\ar@<-3pt>[l]_-{\Phi_{\!\K}} 
} &\!\!\! ,
\end{array}
\end{equation*}

\noindent which can also be deduced from the general 2-adjunction given by \eqref{1405261543}. In this sense, we provide only a few remarks on the structure for the several objects that build this 2-adjunction. \\

The description of $2$-functor, $\Psi_{\K}$, is given completely in order to provide the necessary notation. The structure of such 2-functor goes as follows\\

\begin{enumerate}

\item [1.-] On $0$-cells, $\Psi_{\K}(\mathcal{C}, F) = G_{\F}\dashv V_{\F}$, \emph{i.e.} the Kleisli adjunction.

\item [2.-] On $1$-cells, $(P,\pi) : (\mathcal{C}, F)\longrightarrow  (\mathcal{D}, H)$, $\Psi_{\K}(P, \pi) = (P, P_{\pi}, \rho_{\pi})$. In the  definition of the functor $P_{\pi}:\mathcal{C}_{\F}\longrightarrow \mathcal{D}_{\Hh}$, we use the notation $(\cdot)^{\sharp}$ given for a  morphism in $\mathcal{C}_{\F}$ and $(\cdot)^{\flat}$ for a morphism in  $\mathcal{D}_{\Hh}$. This notation is used in \cite{masa_cawm_2ned} and  \cite{tami_psdl}.

\begin{enumerate}

\item [(i)] On objects, $X$ in $\mathcal{C}_{\F}$, $P_{\pi}X= PX$.

\item [(ii)] On morphisms, $x^{\sharp}:X\longrightarrow Y$ in $\mathcal{C}_{\F}$, $P_{\pi}x^{\sharp}= (\pi C_{x^{\sharp}}\cdot Px)^{\flat}$, where $C_{x^{\sharp}}$ is the notation for the codomain of the morphism $x^{\sharp}$ as in $C_{\F}$, which in this case is $Y$.

\item [(iii)] In order to define $\rho_{\pi}$ we have to prove that the following equality of functors takes place, $G_{\Hh}P = P_{\pi}G_{\F}$.

On objects and morphisms $f:A\longrightarrow B$ in $\mathcal{C}$,

\begin{eqnarray*}
G_{\Hh}PA &=& PA = P_{\pi}A = P_{\pi}G_{\F}A ,\\
G_{\Hh}Pf  &=& (HPf\cdot \eta^{\Hh}PA)^{\flat} = (HPf\cdot \pi A\cdot P\eta^{\F}A)^{\flat} \\
&=& (\pi B\cdot PFf\cdot P\eta^{\F}A)^{\flat} = P_{\pi}(Ff\cdot \eta^{\F}A)^{\sharp} = P_{\pi}G_{\F}f
\end{eqnarray*}

\noindent where the second equality takes place because of the unitality condition on $\pi$ and the third one is due to the naturality on $\pi$. Using \eqref{1407052237}, we get the mate for this identity

\begin{equation*}\Label{1310291824}
\rho_{\pi} = V_{\Hh}P_{\pi}\varepsilon^{DU^{\F}}\scirc \:\eta^{\Hh}PV_{\F},
\end{equation*}

\noindent whose component, at $X$ in $\mathcal{C}_{\F}$, is $\rho_{\pi}X = \mu^{\Hh}PX\cdot H\pi X\cdot \eta^{\Hh}PFX = \pi X$.

\end{enumerate}

\item [3.-] On $2$-cells, $\vartheta:(P,\pi)\longrightarrow(Q,\tau)$, we have

\begin{equation*}\Label{1311011235}
\Psi_{\K}(\vartheta) = (\alpha_{\vartheta}, \beta_{\vartheta})
\end{equation*}

\noindent where $\alpha_{\vartheta}:= \vartheta$ and we rename $\beta_{\vartheta}$ as $\tilde{\vartheta}$. The induced natural tranformation $\tilde{\vartheta}:P_{\pi}\longrightarrow Q_{\tau}:\mathcal{C}_{\F}\longrightarrow \mathcal{D}_{\Hh}$ is defined through its component on $X$, using the condition $G_{\Hh}\separ\vartheta=\tilde{\vartheta}\separ G_{\F}$, as

\begin{equation}\Label{1311011822}
\tilde{\vartheta}X = (\eta^{\Hh}QX\cdot \vartheta X)^{\flat}
\end{equation}

\end{enumerate}

Since we have a 2-adjunction, the following isomorphism of categories takes place, natural in $(\mathcal{X}, H)$ and $L\dashv R$

\begin{equation}\Label{1311042204}
Hom_{\Mnd^{\scriptscriptstyle{\bullet}}({}_{2}Cat)}\big(\: (\mathcal{X}, H), \Phi_{\K}(L\dashv R)\: \big)\cong Hom_{\Adj_{\Ll}({}_{2}Cat)}\big(\:\Psi_{\K}(\mathcal{X}, H),L\dashv R\:\big)
\end{equation}


\section{Monoidal Extensions (Kleisli Type)}

\subsection{Lax Monads}

Dual to colax monads, we give the definition of a \emph{lax monad}.\\

\newtheorem{1403061850}{Definition}[section]
\begin{1403061850}
A lax monad $((F, \zeta, \omega), \mu^{\F}, \eta^{\F})$ over a monoidal category $(\mathcal{C}, \otimes, I, a, l, r)$ consists of the following
\begin{enumerate}
\item $(F,\mu^{\F}, \eta^{\F})$ is a monad on $\mathcal{C}$.
\item $(F,\zeta,\omega):(\mathcal{C}, \otimes, I)\longrightarrow (\mathcal{C}, \otimes, I)$ is a lax monoidal functor. This means that the natural transformations $\zeta:\otimes\cdot (F\times F)\longrightarrow F\cdot\otimes$ and $\omega:\delta_{I}\longrightarrow F\cdot\delta_{I}$, fulfills the commutativity on the following diagrams

\begin{equation}\Label{1403061941}
\xymatrix@C=1.8cm@R=1.4cm{
(FA\otimes FB)\otimes FC\ar[r]^{\zeta_{A,B}\otimes FC}\ar[d]_{a_{FA,FB,FC}} & F(A\otimes B)\otimes FC\ar[r]^{\zeta_{A\otimes B, C}} & F((A\otimes B)\otimes C)\ar[d]^{Fa_{A,B,C}} \\
FA\otimes (FB\otimes FC)\ar[r]_{FA\otimes \zeta_{B,C}} & FA\otimes F(B\otimes C) \ar[r]_{\zeta_{A,B\otimes C}} & F(A\otimes (B\otimes C))
}
\end{equation}

\begin{equation}\Label{1403131301}
\begin{array}{ccc}
\xymatrix@C=1.2cm@R=1.2cm{
I\otimes FA\ar[r]^-{\omega\separ\otimes\separ FA}\ar@/_0.8pc/[rd]_{l_{FA}} & FI \otimes FA\ar[r]^-{\zeta_{I, A}} & F(I\otimes A)\ar@/^0.8pc/[ld]^{Fl_{A}} \\
 & FA  & 
} &
&
\xymatrix@C=1.2cm@R=1.2cm{
F(A\otimes I)\ar@/_0.8pc/[rd]_{Fr_{A}} & FA\otimes FI\ar[l]_-{\zeta_{A,I}} & FA\otimes I\ar[l]_-{FA\separ\otimes\separ \omega}\ar@/^0.8pc/[ld]^{r_{FA}}\\
 & FA & 
}
\end{array}
\end{equation}

\item $\mu^{\F}:(F,\zeta,\omega)\cdot(F,\zeta,\omega)\longrightarrow (F,\zeta,\omega)$ and $\eta^{\F}:(1_{\mathcal{C}}, 1_{\otimes}, 1_{\delta_{I}})\longrightarrow (F,\zeta,\omega)$ are lax natural transformations, the adjective lax adds, to the naturality, the following commutative diagrams

\begin{equation}\Label{1403131208}
\begin{array}{ccc}
\xymatrix@C=1.5cm@R=1.5cm{
\otimes (FF\times FF) \ar[r]^-{\zeta (F\times F)}\ar[d]_{\otimes(\mu^{\F}\times \mu^{\F})} & F\otimes (F\times F)\ar[r]^-{F\zeta } & FF\otimes \ar[d]^{\mu^{\F}\otimes}\\
\otimes (F\times F)\ar[rr]_{\zeta} & &     F\otimes
} & 
\xy<1cm,0cm>
\POS (0,-20.5) *+{,},
\endxy
&
\xymatrix@C=1.15cm@R=1.5cm{
\delta_{I} \ar[r]^-{\omega}\ar@/_1pc/[rrd]_{\omega} & F\delta_{I}\ar[r]^-{F\omega} & FF\delta_{I}\ar[d]^{\mu^{\F}\delta_{I}}\\
 & &     F\delta_{I}
}
\end{array}
\end{equation}

\noindent \&

\begin{equation}\Label{1403131209}
\begin{array}{ccc}
\xy
\POS (-22, 0) *+{\otimes} = "a",
\POS (22, 0) *+{\otimes} = "d",
\POS (22, -20) *+{F\otimes} = "c", 
\POS (-22, -20) *+{\otimes(F\times F)} = "b",
\POS "a" \ar^{1_{\otimes}} "d",
\POS "d" \ar^{\eta^{\F}\otimes} "c",
\POS "a" \ar_{\otimes(\eta^{\F}\times\eta^{\F})} "b",
\POS "b" \ar_{\zeta} "c",
\endxy
& 
\xy<1cm,0cm>
\POS (0,-20.5) *+{,},
\endxy
&
\xy<1cm,0cm>
\POS (0, -4)
\xymatrix@C=1.2cm@R=1cm{
\delta_{I}\ar@/_0.7pc/[rd]_{\omega}\ar[r]^{1_{\delta_{I}}}&  \delta_{I}\ar[d]^{\eta^{\F}\delta_{I}}\\
&     F\delta_{I}
}
\endxy
\end{array}
\end{equation}

\end{enumerate}
\end{1403061850}

\newtheorem{1403131236}[1403061850]{Note}
\begin{1403131236}
Necessarily $\omega(0)= \eta^{\F}\! I$.
\end{1403131236}

The natural transformation $\omega$ has only one component at $0$ in $\mathbf{1}$, then both are going to be denoted by $\omega$, \emph{i.e.} $\omega = \omega(0) = \eta^{\F}\! I$.\\


We are going to make use of the isomorphism \eqref{1311042204}. The result we want to obtain using this isomorphism is the following. 

\newtheorem{1311042045}[1403061850]{Theorem}
\begin{1311042045}

There is a bijective correspondence between the following structures

\begin{enumerate}

\item [1.-] Colax monads $((F, \zeta, \omega),\mu^{\F}, \eta^{\F})$, for the
  monoidal structure $(\mathcal{C}, \otimes, I, a, l, r)$.

\item [2.-] Morphims and transformations of monads of the form

\begin{eqnarray*}
(\otimes, \zeta)&:& (\mathcal{C}\times\mathcal{C}, F\times F)\longrightarrow (\mathcal{C}, F),\\
(\delta_{I}, \omega)&:&(\mathbf{1}\:, 1_{\mathbf{1}}\:)\longrightarrow (\mathcal{C},F)\\
a&:&(\otimes\cdot(\otimes\times
\mathcal{C}),\zeta(\otimes\times\mathcal{C})\circ\otimes(\zeta\times
F))\longrightarrow(\otimes\cdot (\mathcal{C}\times\otimes)\cdot a_{\mathcal{C}}, \zeta(\mathcal{C}\times\otimes)a_{\mathcal{C}}\circ\otimes(F\times\zeta)a_{\mathcal{C}})\\
&&\quad\quad\quad\quad\quad\quad\quad\quad\quad\quad\quad\quad:((\mathcal{C}\times\mathcal{C})\times\mathcal{C}, (F\times F)\times F)\longrightarrow (\mathcal{C}, F), \\
l&:& (\otimes\cdot(\delta_{I}\times \mathcal{C})\cdot l_{\mathcal{C}}^{-1},\zeta\separ(\delta_{I}\times\mathcal{C})\separ l^{-1}_{\mathcal{C}}\circ
\otimes(\omega\times F)\separ l^{-1}_{\mathcal{C}} )\longrightarrow (1_{\mathcal{C}}, 1_{F}):(\mathcal{C}, F)\longrightarrow (\mathcal{C}, F),\\
r&:&(\otimes\cdot(\mathcal{C}\times\delta_{I})\cdot r_{\mathcal{C}}^{-1},
\zeta\separ(\mathcal{C}\times\delta_{I})\separ r^{-1}_{\mathcal{C}}\circ\otimes(F\times\omega)\separ r^{-1}_{\mathcal{C}}
)\longrightarrow (1_{\mathcal{C}}, 1_{F}):(\mathcal{C}, F)\longrightarrow (\mathcal{C}, F).\\ 
\end{eqnarray*}

\item [3.-] Monoidal structures for the Kleisli category $(\mathcal{C}_{\F}, \widetilde{\otimes} , \tilde{I})$ such that the following diagrams of arrows  and surfaces commute

\begin{equation}\Label{1405081433}
\begin{array}{cc}
\xy<1cm,0cm>
\POS (0, 10) *+{(a)},
\POS (0, 0) *+{\mathcal{C}\times\mathcal{C}} = "c1",
\POS (25, 0) *+{\mathcal{C}} = "c2", 
\POS (0, -20) *+{\mathcal{C}_{\F}\times\mathcal{C}_{\F}} = "d1",
\POS (25, -20) *+{\mathcal{C}_{\F}} = "d2", 
\POS (18, -10) = "a",
\POS (10, -18) = "b",
\POS "c1" \ar^-{\otimes} "c2",
\POS "c1" \ar_{G_{\F}\times G_{\F}} "d1",
\POS "c2" \ar^{G_{\F}} "d2",
\POS "d1" \ar_-{\widetilde{\otimes}} "d2",
\endxy & 
\xy<1cm,0cm>
\POS (0, 10) *+{(b)},
\POS (0,0) *+{\mathbf{1}} = "c1",
\POS (25,0) *+{\mathcal{C}} = "c2",
\POS (25,-20) *+{\mathcal{C}^{\F}} = "d2",
\POS (0,-20) *+{\mathbf{1}_{1_{\mathbf{1}}}} = "d1",
\POS "c1" \ar^-{\delta_{I}} "c2",
\POS "c1" \ar_{G_{1_{\mathbf{1}}}} "d1",
\POS "c2" \ar^{G^{\F}} "d2",
\POS "d1" \ar_-{\delta_{\tilde{I}}} "d2",
\endxy
\end{array}
\end{equation}

\begin{equation*}\Label{1405072237}
\begin{array}{ccc}
\xy<1cm,0cm>
\POS (0, 0) *+{\mathcal{C}^{3}} = "c1",
\POS (7, 0) *+{\phantom{(\mathcal{C})}}= "c3",
\POS (30, 0) *+{\mathcal{C}} = "c2", 
\POS (11, 4) = "a",
\POS (11, -4) = "b",
\POS (11,-36) = "c",
\POS (11,-44) = "d",

\POS (0,-40) *+{(\mathcal{C}_{\F})^{3}} = "d1",
\POS (30,-40) *+{\mathcal{C}_{\F}} = "d2",

\POS "c1" \ar@/^1.5pc/^{\otimes\cdot (\otimes\times \mathcal{C})} "c2",
\POS "c1" \ar@/_1.5pc/_{\otimes\cdot (\mathcal{C}\times\otimes)\cdot a_{\mathcal{C}}} "c2",
\POS "a" \ar^{\phantom{11}a} "b",

\POS "d1" \ar@/^1.5pc/^{\widetilde{\otimes}\cdot (\widetilde{\otimes}\times \mathcal{C}_{\F})} "d2",
\POS "d1" \ar@/_1.5pc/_{\widetilde{\otimes}\cdot (\mathcal{C}_{\F}\times\tilde{\otimes})\cdot a_{\mathcal{C}_{\F}}} "d2",
\POS "c" \ar^{\phantom{11}\tilde{a}} "d",

\POS "c1" \ar_{(G_{\F})^{3}} "d1",
\POS "c2" \ar^{G_{\F}} "d2",
\endxy &
\xy<1cm,0cm>
\POS (0, 0) *+{\mathcal{C}} = "c1",
\POS (7, 0) *+{\phantom{(\mathcal{C})}}= "c3",
\POS (30, 0) *+{\mathcal{C}} = "c2", 
\POS (11, 4) = "a",
\POS (11, -4) = "b",
\POS (11,-36) = "c",
\POS (11,-44) = "d",

\POS (0,-40) *+{\mathcal{C}_{\F}} = "d1",
\POS (30,-40) *+{\mathcal{C}_{\F}} = "d2",

\POS "c1" \ar@/^1.5pc/^{\otimes\cdot (\delta_{I}\times \mathcal{C})\cdot l^{-1}_{\mathcal{C}}} "c2",
\POS "c1" \ar@/_1.5pc/_{1_{\mathcal{C}}} "c2",
\POS "a" \ar^{\phantom{12}l} "b",

\POS "d1" \ar@/^1.5pc/^{\widetilde{\otimes}\cdot (\delta_{\tilde{I}}\times
  \mathcal{C}_{\F})\cdot l^{-1}_{\mathcal{C}_{\F}}} "d2",
\POS "d1" \ar@/_1.5pc/_{1_{\mathcal{C}_{\F}}} "d2",
\POS "c" \ar^{\phantom{12}\tilde{l}} "d",

\POS "c1" \ar_{G_{\F}} "d1",
\POS "c2" \ar^{G_{\F}} "d2",
\endxy
&
\xy<1cm,0cm>
\POS (0, 0) *+{\mathcal{C}} = "c1",
\POS (7, 0) *+{\phantom{(\mathcal{C})}}= "c3",
\POS (30, 0) *+{\mathcal{C}} = "c2", 
\POS (11, 4) = "a",
\POS (11, -4) = "b",
\POS (11,-36) = "c",
\POS (11,-44) = "d",

\POS (0,-40) *+{\mathcal{C}_{\F}} = "d1",
\POS (30,-40) *+{\mathcal{C}_{\F}} = "d2",

\POS "c1" \ar@/^1.5pc/^{\otimes\cdot (\mathcal{C}\times\delta_{I})\cdot r^{-1}_{\mathcal{C}}} "c2",
\POS "c1" \ar@/_1.5pc/_{1_{\mathcal{C}}} "c2",
\POS "a" \ar^{\phantom{12}r} "b",

\POS "d1" \ar@/^1.5pc/^{\widetilde{\otimes}\cdot
  (\mathcal{C}_{\F}\times\delta_{\tilde{I}})\cdot r^{-1}_{\mathcal{C}_{\F}}} "d2",
\POS "d1" \ar@/_1.5pc/_{1_{\mathcal{C}_{\F}}} "d2",
\POS "c" \ar^{\phantom{12}\tilde{r}} "d",

\POS "c1" \ar_{G_{\F}} "d1",
\POS "c2" \ar^{G_{\F}} "d2",
\endxy
\end{array}
\end{equation*}

\end{enumerate}

\end{1311042045}

\emph{Proof}:\\

$1\Rightarrow 2$)\\

Consider a lax monad $((F, \zeta, \omega),\mu^{\F}, \eta^{\F})$ for the
monoidal category $(\mathcal{C}, \otimes, I)$. In particular, $\mu^{\F}$ and
$\eta^{\F}$ are natural lax monoidal transformations. Therefore, the
commutativity of the first diagram in \eqref{1403131208} and the first one in
\eqref{1403131209} is equivalent to the condition that the following be a
monad morphism $(\otimes, \zeta):(\mathcal{C}\times\mathcal{C} , F\times F)\longrightarrow (\mathcal{C}, 
F)$.\\

The commutativity condition on the second diagrams in \eqref{1403131208} and
\eqref{1403131209} is equivalent to the condition for the following to be a monad
morphism $(\delta_{I},\omega): (\mathbf{1}\:,
1_{\mathbf{1}}\:)\longrightarrow (\mathcal{C},F) $. \\

Since $(\otimes, \zeta)$ is a morphism of monads so are $(\otimes\cdot
(\otimes\times \mathcal{C}),\:\zeta(\otimes\times
\mathcal{C})\circ\otimes(\zeta\times F))$ and $(\otimes\cdot
(\mathcal{C}\times\otimes)\cdot
a_{\mathcal{C}},\:\zeta(\mathcal{C}\times\otimes)a_{\mathcal{C}}\circ\otimes(F\times\zeta)a_{\mathcal{C}})$. Yet
again, since 
$((F, \zeta), \mu^{\F}, \eta^{\F})$ is a lax monad over the monoidal category
$(\mathcal{C},\otimes, I, a, l, r)$, then a commutative diagram like
\eqref{1403061941} takes place. Therefore the following is a 2-cell in $\Mnd^{\scriptscriptstyle{\bullet}}({}_{2}Cat)$.

\begin{equation*}
\xy<1cm,0cm>
\POS (0, 0) *+{((\mathcal{C}\times\mathcal{C})\times\mathcal{C}, (F\times F)\times F)} = "c1",
\POS (7, 0) *+{\phantom{(\mathcal{C})}}= "c3",
\POS (50, 0) *+{(\mathcal{C}, F)} = "c2", 
\POS (27, 5) = "a",
\POS (27, -5) = "b",
\POS "c3" \ar@/^2pc/^{(\otimes\cdot (\otimes\times \mathcal{C}),\:\zeta(\otimes\times \mathcal{C})\circ\otimes(\zeta\times F))} "c2",
\POS "c3" \ar@/_2pc/_{(\otimes\cdot (\mathcal{C}\times\otimes)\cdot a_{\mathcal{C}},\:\zeta(\mathcal{C}\times\otimes)a_{\mathcal{C}}\circ\otimes(F\times\zeta)a_{\mathcal{C}})} "c2",
\POS "a" \ar^{\phantom{1}a} "b",
\endxy
\end{equation*}\\

Since $(\otimes, \zeta)$ and $(\delta_{I}, \omega)$ are monad morphisms so is $(\otimes\cdot(\delta_{I}\times \mathcal{C})\cdot l_{\mathcal{C}}^{-1},\zeta\separ(\delta_{I}\times\mathcal{C})\separ l^{-1}_{\mathcal{C}}\circ
\otimes\separ(\omega\times F)\separ l^{-1}_{\mathcal{C}})$ and taking into
account the commutativity of the diagram  (\ref{1403131301}a), we can state
that the following is a 2-cell in $\Mnd^{\scriptscriptstyle{\bullet}}({}_{2}Cat)$

\begin{equation*}
\xy<1cm,0cm>
\POS (0, 0) *+{(\mathcal{C}, F)} = "c1",
\POS (40, 0) *+{(\mathcal{C}, F)} = "c2", 
\POS (13, 3.5) = "a",
\POS (13, -3.5) = "b",
\POS "c1" \ar@/^1.5pc/^{(\otimes\cdot (\delta_{I}\times\mathcal{C})\cdot l^{-1}_{\mathcal{C}}, \:\zeta\separ(\delta_{I}\times\mathcal{C})\separ l^{-1}_{\mathcal{C}}\circ\separ
\otimes\separ(\omega\times F)\separ l^{-1}_{\mathcal{C}})} "c2",
\POS "c1" \ar@/_1.5pc/_{(1_{\mathcal{C}},1_{F})} "c2",
\POS "a" \ar^{\phantom{111}l} "b",
\endxy
\end{equation*}

In the very same way, the following is a $2$-cell of monads,
$r:(\otimes\cdot(\mathcal{C}\times\delta_{I})\cdot
r^{-1}_{\mathcal{C}},\zeta(\mathcal{C}\times\delta_{I})\separ
r^{-1}_{\mathcal{C}}\circ\otimes(F\times\omega)\separ r^{-1}_{\mathcal{C}}
)$\\ 

2 $\Rightarrow$ 1) The previous assertions can be reverted.\\

2 $\Rightarrow$ 3)\\ 

Suppose we have a monad morphism $(\otimes, \zeta)$. Use the isomorphism \eqref{1311042204}, with $(\mathcal{D}, H,\mu^{\Hh},\eta^{\Hh}) = (\mathcal{C}\times\mathcal{C}, F\times F, \mu^{\F}\times \mu^{\F}, \eta^{\F}\times \eta^{\F})$ and $L\dashv R = G_{\F}\dashv V_{\F}$ to get an associated morphism of adjunctions $(\otimes, \otimes_{\zeta}): G_{\F}\times G_{\F}\dashv V_{\F}\times V_{\F}\longrightarrow G_{\F}\dashv V_{\F}$, such that a diagram like (\ref{1405081433}a) commutes. According to the definition of $\Psi_{\K}$, the functor $\otimes_{\zeta}$ acts as follows. On objects, 

\begin{equation*}\Label{1403131226}
\otimes_{\zeta}(X,Y) = \otimes(X, Y) = X\otimes Y,
\end{equation*}

\noindent and on morphisms,

\begin{equation*}\Label{1403131227}
\otimes_{\zeta}(x^{\sharp},y^{\sharp}) = (\zeta_{C_{x^{\sharp}},C_{y^{\sharp}}}\cdot (x\otimes y))^{\sharp}
\end{equation*}\\

\noindent where $C_{x^{\sharp}}$ is codomain of the morphism $x^{\sharp}$ for example. We rename $\otimes_{\zeta}$ as $\widetilde{\otimes}$.\\

For the monad morphism, $(\delta_{I},\omega):(\mathbf{1},1_{\mathbf{1}})\longrightarrow (\mathcal{C}, F)$, use the mentioned
isomorphism with $(\mathcal{D}, H,\mu^{\Hh},\eta^{\Hh}) = (\mathbf{1}\:,1_{\mathbf{1}}, 1_{1_{\mathbf{1}}}, 1_{1_{\mathbf{1}}})$,
\emph{i.e.} the trivial monad on the category $\mathbf{1}$, and $L\dashv R = G_{\F}\dashv V_{\F}$. Therefore, there exists
an adjunction morphism $(\delta_{I}, \lbrack\delta_{I}\rbrack_{\omega}): G_{1_{\mathbf{1}}}\dashv V_{1_{\mathbf{1}}}\longrightarrow G_{\F}\dashv V_{\F}$. According to the 2-functor $\Psi_{\K}$, the functor $\lbrack\delta_{I}\rbrack_{\omega}:\mathbf{1}\longrightarrow \mathcal{C}_{\F}$, acts in the following way

\begin{equation*}\Label{1403131233}
\lbrack\delta_{I}\rbrack_{\omega}(0) = \delta_{I}(0) = I
\end{equation*}

\begin{equation*}\Label{1403131234}
\lbrack\delta_{I}\rbrack_{\omega}(1_{0}) = (\omega(0)\cdot \delta_{I}(1_{0}))^{\sharp} = (\eta^{\F}I)^{\sharp}
\end{equation*}\\

That is to say $\lbrack\delta_{I}\rbrack_{\omega} = \delta_{\tilde{I}}:\mathbf{1}_{1_{\mathbf{1}}}\longrightarrow \mathcal{C}_{\F}$, where $\tilde{I} = I$.\\

Suppose that we have the following $2$-cell in $\Mnd^{\bullet}({}_{2}Cat)$,\\

\begin{equation*}
\xy<1cm,0cm>
\POS (0, 0) *+{((\mathcal{C}\times\mathcal{C})\times\mathcal{C}, (F\times F)\times F)} = "c1",
\POS (7, 0) *+{\phantom{(\mathcal{C})}}= "c3",
\POS (50, 0) *+{(\mathcal{C}, F)} = "c2", 
\POS (27, 5) = "a",
\POS (27, -5) = "b",
\POS "c3" \ar@/^2pc/^{(\otimes\cdot (\otimes\times \mathcal{C}),\:\zeta(\otimes\times \mathcal{C})\circ\otimes(\zeta\times F))} "c2",
\POS "c3" \ar@/_2pc/_{(\otimes\cdot (\mathcal{C}\times\otimes)\cdot
  a_{\mathcal{C}},
  \:\zeta(\mathcal{C}\times\otimes)a_{\mathcal{C}}\circ\otimes(F\times\zeta) a_{\mathcal{C}})} "c2",
\POS "a" \ar^{\phantom{1}a} "b",
\endxy
\end{equation*}\\

In order to continue with the calculations, we use the following notation, for
the  sake of simplification\\

\begin{eqnarray*}
{}^{\cdot}\zeta^{2} &:=& \zeta(\otimes\times \mathcal{C})\circ\otimes(\zeta\times F),\\
{}_{\cdot}\zeta^{2} &:=& \zeta(\mathcal{C}\times\otimes)a_{\mathcal{C}}\circ\otimes(F\times\zeta)a_{\mathcal{C}}, \\
(\cdot)^{3}               &:=& (\cdot\times\cdot)\times\cdot\:.
\end{eqnarray*}

According to the isomorphism of categories given by \eqref{1311042204}, to the previous $2$-cell in $\Mnd^{\scriptscriptstyle{\bullet}}({}_{2}Cat)$ corresponds a $2$-cell, $(\alpha_{a}, \beta_{a})$ in $\Adj_{\Ll}({}_{2}Cat)$, where $\alpha_{a} = a$ and we rename $\beta_{a}= \tilde{a}$ and such that 

\begin{equation*}
\xy<1cm,0cm>
\POS (0, 0) *+{\mathcal{C}^{3}} = "c1",
\POS (7, 0) *+{\phantom{(\mathcal{C})}}= "c3",
\POS (40, 0) *+{\mathcal{C}} = "c2", 
\POS (14, 4) = "a",
\POS (14, -4) = "b",
\POS (14,-36) = "c",
\POS (14,-44) = "d",

\POS (0,-40) *+{(\mathcal{C}_{\F})^{3}} = "d1",
\POS (40,-40) *+{\mathcal{C}_{\F}} = "d2",

\POS "c1" \ar@/^1.7pc/^{\otimes\cdot (\otimes\times \mathcal{C})} "c2",
\POS "c1" \ar@/_1.7pc/_{\otimes\cdot (\mathcal{C}\times\otimes)\cdot a_{\mathcal{C}}} "c2",
\POS "a" \ar^{\phantom{12}a} "b",

\POS "d1" \ar@/^1.7pc/^{[\otimes\cdot (\otimes\times \mathcal{C})]_{{}^{\cdot}\!\zeta^{2}}} "d2",
\POS "d1" \ar@/_1.7pc/_{[\otimes\cdot (\mathcal{C}\times\otimes)\cdot a_{\mathcal{C}}]_{{}_{\cdot}\zeta^{2}}} "d2",
\POS "c" \ar^{\phantom{12}\widetilde{a}} "d",

\POS "c1" \ar@<-4pt>_{(G_{\F})^{3}} "d1",
\POS "d1" \ar@<-1pt>_{(V_{\F})^{3}} "c1",
\POS "c2" \ar@<-2pt>_{G_{\F}} "d2",
\POS "d2" \ar@<-3pt>_{V_{\F}} "c2",
\endxy
\end{equation*}\\

It can be show that 

\begin{eqnarray*}
\lbrack\otimes\cdot (\otimes\times \mathcal{C})\rbrack_{{}^{\cdot}\zeta^{2}} &=& \widetilde{\otimes}\cdot(\widetilde{\otimes}\times \mathcal{C}_{\F})\\
\lbrack\otimes\cdot (\mathcal{C}\times\otimes)\cdot a_{\mathcal{C}} \rbrack_{{}_{\cdot}\zeta^{2}} &=& \widetilde{\otimes}\cdot (\mathcal{C}_{\F}\times\widetilde{\otimes})\cdot a_{\mathcal{C}_{\F}}
\end{eqnarray*}\\

Therefore, we have a natural transformation $\widetilde{a}:\widetilde{\otimes}\cdot(\widetilde{\otimes}\times\mathcal{C}_{\F})\longrightarrow \widetilde{\otimes}\cdot(\mathcal{C}_{\F}\times\widetilde{\otimes})\cdot a_{\mathcal{C}_{\F}}$ that will
be part of a monoidal structure  on
$\mathcal{C}_{\F}$. According to the 2-functor $\Psi_{\K}$, the component of $\widetilde{a}$ at $((X, Y), Z)$ is\\

\begin{equation*}\Label{1405091636}
\widetilde{a}_{X,Y,Z} = (\eta^{\F}(X\otimes (Y\otimes Z))\cdot a_{X,Y,Z})^{\sharp} 
\end{equation*}\\

Suppose that we have a 2-cell in $\Mnd^{\scriptscriptstyle{\bullet}}({}_{2}Cat)$ of the form
$l:\big(\otimes\cdot(\delta_{I}\times\mathcal{C})\cdot l^{-1}_{\mathcal{C}}, \zeta(\delta_{I}\times\mathcal{C})\separ l^{-1}_{\mathcal{C}}\circ\otimes\separ(\omega\times F)\separ l^{-1}_{\mathcal{C}}\big)\longrightarrow (1_{\mathcal{C}}, 1_{\F}):(\mathcal{C}, F)\longrightarrow (\mathcal{C}, F)$. \\

Therefore, we obtain a natural transformation $\tilde{l}:\widetilde{\otimes}\cdot (\delta_{\tilde{I}}\times\mathcal{C}_{\F})\cdot
l^{-1}_{\mathcal{C}_{\F}}\longrightarrow 1_{\mathcal{C}_{\F}}$. Using the
definition of the functor $\Psi_{\K}$ on the 2-cell $l$, the component of $\tilde{l}$, on the object $X$ in $\mathcal{C}_{\F}$, is \\

\begin{equation}\label{1405091630}
\tilde{l}X = (\eta^{\F}X\cdot lX)^{\sharp}
\end{equation}\\

Similarly, for the monad morphism
$r:(\otimes\cdot(\mathcal{C}\times\delta_{I})\cdot r^{-1}_{\mathcal{C}}, \zeta
(\mathcal{C}\times\delta_{I})r^{-1}_{\mathcal{C}}\circ\otimes(F\times\omega)r^{-1}_{\mathcal{C}})\longrightarrow
(1_{\mathcal{C}}, 1_{\F}):(\mathcal{C}, F)\longrightarrow(\mathcal{C}, F)$, we
obtain a natural transformation $\tilde{r}:
\otimes_{\zeta}\cdot(\mathcal{C}_{\F}\times\delta_{\tilde{I}})\cdot
r^{-1}_{\mathcal{C}_{\F}}\longrightarrow 1_{\mathcal{C}_{\F}}:
\mathcal{C}_{\F}\longrightarrow \mathcal{C}_{\F}$.\\

The proof of the coherence conditions are left to the reader.\\

In summary, $(\mathcal{C}_{\F}, \widetilde{\otimes}, \tilde{I}, \tilde{a}, \tilde{l}, \tilde{r})$ has a monoidal structure on $\mathcal{C}_{\F}$.\\

3 $\Rightarrow$ 2)\\

Using the isomorphism, given by \eqref{1311042204}, we get the return of the 
proof. For example, the image, under $\Phi_{\K}$, for the 2-cell of adjunctions $(a,\tilde{a}):(\otimes\cdot(\otimes\times\mathcal{C}),
\otimes\cdot(\mathcal{C}\times\otimes)\cdot a_{\mathcal{C}})\longrightarrow(\widetilde{\otimes}\cdot(\widetilde{\otimes}\times\mathcal{C}_{\F}),\widetilde{\otimes}\cdot(\mathcal{C}_{\F}\times\widetilde{\otimes})\cdot a_{\mathcal{C}_{\F}}): (G_{\F}\times G_{\F})\times G_{\F}\dashv (V_{\F}\times V_{\F})\times V_{\F}$ is 

\begin{eqnarray*}
\Phi_{\K}((a,\tilde{a})) &=& a:(\otimes,
\pi_{\otimes})(\otimes\times\mathcal{C},\pi_{\otimes\times\mathcal{C}})\longrightarrow
(\otimes, \pi_{\otimes})(\mathcal{C}\times\otimes,
\pi_{\mathcal{C}\times\otimes})(a_{\mathcal{C}},
\pi_{a_{\mathcal{C}}}):(\mathcal{C}^{3}, F^{3})\longrightarrow (\mathcal{C},
F)\\
&=& a:(\otimes,
\pi_{\otimes})(\otimes\times\mathcal{C},\pi_{\otimes}\times
F)\longrightarrow (\otimes,
\pi_{\otimes})(\mathcal{C}\times\otimes,F\times\pi_{\otimes})(a_{\mathcal{C}},
1_{\F\times(\F\times\F)\cdot a_{\mathcal{C}}}):\\
&& (\mathcal{C}^{3},
F^{3})\longrightarrow (\mathcal{C}, F)\\
&=&
a:\big(\otimes\cdot(\otimes\times\mathcal{C}), \pi_{\otimes}(\otimes\times\mathcal{C})\circ\otimes(\pi_{\otimes}\times
F)\big)\longrightarrow\\
&& \big(\otimes\cdot(\mathcal{C}\times\mathcal{C})\cdot a_{\mathcal{C}},\pi_{\otimes}(\mathcal{C}\times\otimes)a_{\mathcal{C}}\circ
\otimes(F\times\pi_{\otimes})a_{\mathcal{C}}\big):(\mathcal{C}^{3},
F^{3})\longrightarrow (\mathcal{C}, F)
\end{eqnarray*}\\

We used the fact that $a_{\mathcal{C}}$ is a morphism of adjunctions. 

\begin{flushright}
$\square$
\end{flushright}

\section{Liftings to the Eilenberg-Moore algebras \& Extensions to the Kleisli Categories}

This is probably the most explored section in this article, a few examples
of the detailed proofs for the following statements are found in
\cite{bofr_haca_II} and \cite{tami_psdl}. In this section, we treated these statements only as direct
consequences of the isomorphisms of categories given by \eqref{1403191902} and \eqref{1311042204}.

\newtheorem{1405102218}{Theorem}[section]
\begin{1405102218}
There is a bijective correspondence, natural in $(\mathcal{C}, F, \mu^{\F}, \eta^{\F})$ and $(\mathcal{D}, H, \mu^{\Hh}, \eta^{\Hh})$,  between the following structures

\begin{enumerate}

\item [1.-] Liftings to the Eilenberg-Moore algebras, for the functor
  $P:\mathcal{C}\longrightarrow \mathcal{D}$. That is to say, the following diagram commutes

\begin{equation*}
\xy<1cm,0cm>
\POS (0, 0) *+{\mathcal{C}} = "a21",
\POS (20, 0) *+{\mathcal{D}}= "a22",
\POS (0, 20) *+{\mathcal{C}^{\F}} = "a11",
\POS (20,20) *+{\mathcal{D}^{\Hh}} = "a12",

\POS "a11" \ar^{Q} "a12",
\POS "a21" \ar_{P} "a22",
\POS "a11" \ar_{U^{\F}} "a21",
\POS "a12" \ar^{U^{\Hh}} "a22",

\endxy
\end{equation*}

\item [2.-] Morphisms of monads $(P,\varphi) : (\mathcal{C}, F)\longrightarrow (\mathcal{D}, H)$. That is to say, a natural transformation  $\varphi:HP\longrightarrow PF$, such that the following diagrams commute

\begin{equation*}
\begin{array}{cc}
\xy<1cm,0cm>
\POS (0,22) *+{HHP} = "a11",
\POS (30,22) *+{HPF} = "a12",
\POS (60,22) *+{PFF} = "a13",
\POS (0,0) *+{HP} = "a21",
\POS (60,0) *+{PF} = "a23",

\POS "a11" \ar^{H\varphi} "a12",
\POS "a12" \ar^{\varphi F} "a13",
\POS "a21" \ar_{\varphi} "a23",
\POS "a11" \ar_{\mu^{\Hh}\! P} "a21",
\POS "a13" \ar^{P\mu^{\F}} "a23"
\endxy & 
\xy<1cm,0cm>
\POS (20,15) *+{P} = "a12",
\POS (0,0) *+{HP} = "a21",
\POS (40,0) *+{PF} = "a23",

\POS "a12" \ar_{\eta^{\Hh}\! P} "a21",
\POS "a12" \ar^{P\eta^{\F}} "a23",
\POS "a21" \ar_{\varphi} "a23",
\endxy
\end{array}
\end{equation*}\\

\end{enumerate}

\end{1405102218} 

\newtheorem{1405102219}[1405102218]{Theorem}
\begin{1405102219}
There exists a bijective correspondence, natural in $(\mathcal{C}, F, \mu^{\F}, \eta^{\F})$ and $(\mathcal{D}, H, \mu^{\Hh}, \eta^{\Hh})$, between the following structures

\begin{enumerate}

\item [1.-] Extensions to the Kleisli categories, for the functor  $P:\mathcal{C}\longrightarrow \mathcal{D}$. That is to say, the following diagram commutes 

\begin{equation*}
\xy<1cm,0cm>
\POS (0, 0) *+{\mathcal{C}_{\F}} = "a21",
\POS (20, 0) *+{\mathcal{D}_{\Hh}}= "a22",
\POS (0, 20) *+{\mathcal{C}} = "a11",
\POS (20,20) *+{\mathcal{D}} = "a12",

\POS "a11" \ar^{P} "a12",
\POS "a21" \ar_{Q} "a22",
\POS "a11" \ar_{G_{\F}} "a21",
\POS "a12" \ar^{G_{\Hh}} "a22",

\endxy
\end{equation*}

\item [2.-] Morphisms of monads $(P, \varphi):(\mathcal{C}, F)\longrightarrow  (\mathcal{D}, H)$. That is to say, a natural transformation  $\varphi:PF\longrightarrow HP$

\begin{equation*}
\begin{array}{cc}
\xy<1cm,0cm>
\POS (0,22) *+{PFF} = "a11",
\POS (30,22) *+{HPF} = "a12",
\POS (60,22) *+{HHP} = "a13",
\POS (0,0) *+{PF} = "a21",
\POS (60,0) *+{HP} = "a23",

\POS "a11" \ar^{\varphi F} "a12",
\POS "a12" \ar^{H\varphi } "a13",
\POS "a21" \ar_{\varphi} "a23",
\POS "a11" \ar_{P\mu^{\F}} "a21",
\POS "a13" \ar^{\mu^{\Hh}P} "a23"
\endxy & 
\xy<1cm,0cm>
\POS (20,15) *+{P} = "a12",
\POS (0,0) *+{PF} = "a21",
\POS (40,0) *+{HP} = "a23",

\POS "a12" \ar_{P\eta^{\F}} "a21",
\POS "a12" \ar^{\eta^{\Hh}\! P} "a23",
\POS "a21" \ar_{\varphi} "a23",
\endxy
\end{array}
\end{equation*}\\

\end{enumerate}

\end{1405102219}

\section{Actions on the Kleisli Category}

\subsection{Categorical Actions}

In this section we give the definition of a \emph{categorical action}.

\newtheorem{1407221804}{Definition}[section]
\begin{1407221804}
Let $(\mathcal{C}, \otimes, I)$ be a monoidal category. A left
$\mathcal{C}$-action on the category $\mathcal{B}$ is a functor
$\boxtimes:\mathcal{C}\times\mathcal{B}\longrightarrow \mathcal{B}$ together with natural
transformations
$\nu:\boxtimes(\otimes\times\mathcal{B})\longrightarrow\boxtimes(\mathcal{C}\times\otimes)a_{\ast}:
(\mathcal{C}\times\mathcal{C})\times \mathcal{B}\longrightarrow \mathcal{B}$
and $j:\boxtimes(\delta_{I}\times\mathcal{B})l^{-1}_{\mC}\longrightarrow
1_{\mathcal{B}}:\mathcal{B}\longrightarrow \mathcal{B}$ such that they fulfill the following commutative diagrams, for
objects $C,C', C''$ in $\mathcal{C}$ and $B$ in $\mathcal{B}$,

\begin{equation*}
\xy<1cm,0cm>
\POS (0, 40) *+{\lbrack (C\otimes C')\otimes C''\rbrack\boxtimes B} = "a11",
\POS (70, 40) *+{(C\otimes C')\boxtimes(C''\boxtimes B)} = "a12",
\POS (0,20) *+{\lbrack C\otimes (C'\otimes C'')\rbrack\boxtimes B} = "a21",
\POS (0,0) *+{C\boxtimes\lbrack (C'\otimes C'')\boxtimes B\rbrack} = "a31",
\POS (70,0) *+{C\boxtimes\lbrack C'\boxtimes (C''\boxtimes B)\rbrack} = "a32",

\POS "a11" \ar^{\nu_{C\otimes C', C'', B}} "a12",
\POS "a12" \ar^{\nu_{C, C', C''\boxtimes B}} "a32",
\POS "a11" \ar_{a_{C, C', C''}\boxtimes B} "a21",
\POS "a21" \ar_{\nu_{C, C'\otimes C'', B}} "a31",
\POS "a31" \ar_{C\boxtimes\nu_{C', C'', B}} "a32",
\endxy
\end{equation*}

\noindent and 

\begin{equation*}
\begin{array}{ccc}
\xy<1cm,0cm>
\POS (0,20) *+{(C\otimes I)\boxtimes B} = "a11",
\POS (40,20)*+{C\boxtimes(I\boxtimes B)} = "a13",
\POS (20, 0)*+{C\boxtimes B} = "a22",

\POS "a11" \ar^{\nu_{C, I, B}} "a13",
\POS "a13" \ar^{C\separ\boxtimes j_{B}} "a22",
\POS "a11" \ar_{r_{C}\boxtimes B} "a22",
\endxy & & 
\xy<1cm,0cm>
\POS (0,20) *+{(I\otimes C)\boxtimes B} = "a11",
\POS (40,20)*+{I\boxtimes(C\boxtimes B)} = "a13",
\POS (20, 0)*+{C\boxtimes B} = "a22",

\POS "a11" \ar^{\nu_{I, C, B}} "a13",
\POS "a13" \ar^{j_{C\separ\boxtimes B}} "a22",
\POS "a11" \ar_{l_{C}\boxtimes B} "a22",
\endxy
\end{array}
\end{equation*}

\end{1407221804}


\subsection{Strong Monads}

In this section we give the definition of a \emph{strong monad}. 

\newtheorem{1405191543}{Definition}[section]
\begin{1405191543}
A \emph{right strong monad} $((F,\sigma^{r}), \mu^{\F}, \eta^{\F})$, on the monoidal category $(\mathcal{C}, \otimes, I)$, is a usual monad $(F, \mu^{\F}, \eta^{\F})$, on $\mathcal{C}$, with a natural transformation $\sigma^{r}:A \otimes FB\longrightarrow F(A\otimes B)$ such that the following diagrams commute

\begin{equation}\Label{1405191712}
\begin{array}{cc}
\xy<1cm,0cm>
\POS (0,20) *+{A\otimes FFB} = "a11",
\POS (35,20) *+{F(A\otimes FB)} = "a12",
\POS (70,20) *+{FF(A\otimes B)} = "a13"
\POS (0,0) *+{A\otimes FB} = "a21",
\POS (70,0) *+{F(A\otimes B)} = "a23",

\POS (-12,25) *+{(a)},

\POS "a11" \ar^{\sigma^{r}_{A,FB}} "a12",
\POS "a12" \ar^{F\sigma^{r}_{\! A,B}} "a13",
\POS "a13" \ar^{\mu^{F}(A\otimes B)} "a23"
\POS "a11" \ar_{A\otimes\mu^{\F}\! B} "a21",
\POS "a21" \ar_{\sigma^{r}} "a23",

\endxy & 
\xy<1cm,0cm>

\POS (20,20) *+{A\otimes B} = "a12",
\POS (0,0) *+{A\otimes FB} = "a21",
\POS (40,0) *+{F(A\otimes B)} = "a23",

\POS (-5,25) *+{(b)},

\POS "a12" \ar^{\eta^{\F}(A\otimes B)} "a23",
\POS "a12" \ar_{A\otimes\eta^{\F}\! B} "a21",
\POS "a21" \ar_{\sigma^{r}} "a23",

\endxy
\end{array}
\end{equation}

\noindent and

\begin{equation}\Label{1405191723}
\begin{array}{cc}
\xy<1cm,0cm>
\POS (0,20) *+{(A\otimes B)\otimes FC} = "a11",
\POS (90,20) *+{F((A\otimes B)\otimes C)} = "a13",
\POS (0,0) *+{A\otimes (B\otimes FC)} = "a21",
\POS (45,0) *+{A\otimes F(B\otimes C)} = "a22",
\POS (90,0) *+{F(A\otimes (B\otimes C))} = "a23"

\POS (0,30) *+{(a)},

\POS "a11" \ar^{\sigma^{r}_{A\otimes B, C}} "a13",
\POS "a13" \ar^{Fa_{A,B,C}} "a23",
\POS "a11" \ar_{a_{A,B,FC}} "a21",
\POS "a21" \ar_{A\otimes\sigma^{r}_{B, C}} "a22",
\POS "a22" \ar_{\sigma^{r}_{A,B\otimes C}} "a23",

\endxy & 
\xy<1cm,0cm>

\POS (0, 20) *+{I\otimes FA} = "a11",
\POS (40,20) *+{F(I\otimes A)} = "a13",
\POS (20,0) *+{FA} = "a22",

\POS (0,30) *+{(b)},

\POS "a11" \ar^{\sigma^{r}_{I, A}} "a13",
\POS "a13" \ar^{Fl_{A}} "a22",
\POS "a11" \ar_{l_{FA}} "a22",

\endxy
\end{array}
\end{equation}

\end{1405191543}

\newtheorem{1405201903}[1405191543]{Definition}
\begin{1405201903}

A \emph{left strong monad} $((F,\sigma^{r}), \mu^{\F}, \eta^{\F})$ on a monoidal category $(\mathcal{C}, \otimes, I)$, is a usual monad $(F, \mu^{\F},\eta^{\F})$ on $\mathcal{C}$, together with a natural transformation $\sigma^{l}_{A,B}:FA\otimes B \longrightarrow F(A\otimes B)$ such that fulfills the commutativity of dual diagrams like \eqref{1405191712} and \eqref{1405191723}.
\end{1405201903}

The following theorem can be stated, note that the proof is just an adaptation
for the corresponding lax monoidal case.

\newtheorem{1505191350}[1405191543]{Theorem}
\begin{1505191350}

There exists a bijection between the following structures 

\begin{enumerate}

\item [1.-] Right strong monads $((F,\sigma^{r}),\mu^{\F},\eta^{\F})$ on the
  monoidal category $(\mathcal{C}, \otimes, I, a, r, l)$.

\item [2.-] Morphisms and transformations of monads of the form

\begin{eqnarray*}
(\otimes, \sigma^{r}) &:& (\mathcal{C}\times\mathcal{C},
\mathcal{C}\times F)\longrightarrow (\mathcal{C}, F)\\
a &:& (\otimes\cdot(\otimes\times\mathcal{C}),
\sigma^{r}(\otimes\times\mathcal{C}))\longrightarrow
(\otimes\cdot(\mathcal{C}\times\otimes)\cdot a_{\mathcal{C}},
\sigma^{r}(\mathcal{C}\times\otimes)a_{\mathcal{C}}\circ
\otimes(\mathcal{C}\times\sigma^{r})a_{\mathcal{C}})\\
&:& ((\mathcal{C}\times\mathcal{C})\times\mathcal{C},
(\mathcal{C}\times\mathcal{C})\times F)\longrightarrow (\mathcal{C}, F)\\
l&:& (\otimes\cdot(\delta_{I}\times\mathcal{C})\cdot l^{-1}_{\mathcal{C}},
\separ\sigma^{r}(\delta_{I}\times\mathcal{C})\separ
l^{-1}_{\mathcal{C}})\longrightarrow (1_{\mathcal{C}}, 1_{\F}):(\mathcal{C},
F)\longrightarrow (\mathcal{C}, F)
\end{eqnarray*}

\item [3.-] Left actions on the Kleisli category, $\mathcal{C}_{\F}$,  $\boxtimes:\mathcal{C}\times\mathcal{C}_{\F}\longrightarrow \mathcal{C}_{\F}$  such that the following diagrams of morphisms and surfaces commute 

\begin{equation}\Label{1405222031}
\xy<1cm,0cm>
\POS (0,20) *+{\mathcal{C}\times\mathcal{C}} = "a11",
\POS (20,20) *+{\mathcal{C}} = "a12",
\POS (0,0) *+{\mathcal{C}\times\mathcal{C}_{\F}} = "a21",
\POS (20,0) *+{\mathcal{C}_{\F}} = "a22",

\POS "a11" \ar^-{\otimes} "a12",
\POS "a12" \ar^{G_{\F}} "a22",
\POS "a11" \ar_{\mathcal{C}\times G_{\F}} "a21",
\POS "a21" \ar_-{\boxtimes} "a22",

\endxy
\end{equation}

\begin{equation}\Label{1405202049}
\begin{array}{ccc}
\xy<1cm,0cm>
\POS (0,15) *+{(a)}
\POS (0, 0) *+{\mathcal{C}^{2}\times\mathcal{C}} = "c1",
\POS (7, 0) *+{\phantom{(\mathcal{C})}}= "c3",
\POS (30, 0) *+{\mathcal{C}} = "c2", 
\POS (11, 4) = "a",
\POS (11, -4) = "b",
\POS (11,-36) = "c",
\POS (11,-44) = "d",

\POS (0,-40) *+{\mathcal{C}^{2}\times\mathcal{C}_{\F}} = "d1",
\POS (30,-40) *+{\mathcal{C}_{\F}} = "d2",

\POS "c1" \ar@/^1.5pc/^{\otimes\cdot (\otimes\times \mathcal{C})} "c2",
\POS "c1" \ar@/_1.5pc/_{\otimes\cdot (\mathcal{C}\times\otimes)\cdot a_{\mathcal{C}}} "c2",
\POS "a" \ar^{\phantom{11}a} "b",

\POS "d1" \ar@/^1.5pc/^{\boxtimes\cdot (\otimes\times \mathcal{C}_{\F})} "d2",
\POS "d1" \ar@/_1.5pc/_{\boxtimes\cdot (\mathcal{C}\times\boxtimes)\cdot a_{\mathcal{C}_{\ast}}} "d2",
\POS "c" \ar^{\phantom{11}\tilde{a}} "d",

\POS "c1" \ar_{\mathcal{C}^{2}\times G_{\F}} "d1",
\POS "c2" \ar^{G_{\F}} "d2",
\endxy & 
\xy<1cm,0cm>
\POS (0, 0) *+{\phantom{123}}
\endxy
&
\xy<1cm,0cm>
\POS (0,15) *+{(b)}
\POS (0, 0) *+{\mathcal{C}} = "c1",
\POS (7, 0) *+{\phantom{(\mathcal{C})}}= "c3",
\POS (30, 0) *+{\mathcal{C}} = "c2", 
\POS (11, 4) = "a",
\POS (11, -4) = "b",
\POS (11,-36) = "c",
\POS (11,-44) = "d",

\POS (0,-40) *+{\mathcal{C}_{\F}} = "d1",
\POS (30,-40) *+{\mathcal{C}_{\F}} = "d2",

\POS "c1" \ar@/^1.5pc/^{\otimes\cdot (\delta_{I}\times\mathcal{C})\cdot l^{-1}_{\mathcal{C}}} "c2",
\POS "c1" \ar@/_1.5pc/_{1_{\mathcal{C}}} "c2",
\POS "a" \ar^{\phantom{12}l} "b",

\POS "d1"
\ar@/^1.5pc/^{\boxtimes\cdot(\delta_{I}\times\mathcal{C}_{\F})\cdot l^{-1}_{\mathcal{C}_{\F}}} "d2",
\POS "d1" \ar@/_1.5pc/_{1_{\mathcal{C}_{\F}}} "d2",
\POS "c" \ar^{\phantom{12}\tilde{l}} "d",

\POS "c1" \ar_{G_{\F}} "d1",
\POS "c2" \ar^{G_{\F}} "d2",
\endxy
\end{array}
\end{equation}

\end{enumerate}

\end{1505191350}

\begin{flushright}
$\square$
\end{flushright}

We state the dual theorem\\

\newtheorem{1505241921}[1405191543]{Theorem}
\begin{1505241921}

There exists a bijection between the following structures 

\begin{enumerate}

\item [1.-] Left strong monads $((F,\sigma^{l}),\mu^{\F},\eta^{\F})$ on the
  monoidal category $(\mathcal{C}, \otimes, I, a, r, l)$.

\item [2.-] Morphisms and transformations of monads of the form

\begin{eqnarray*}
(\otimes, \varphi) &:& (\mathcal{C}\times\mathcal{C},
F\times\mathcal{C})\longrightarrow (\mathcal{C}, F)\\
a &:& (\otimes\cdot(\otimes\times\mathcal{C}),
\varphi(\otimes\times\mathcal{C})\circ\otimes(\varphi\times\mathcal{C}))\longrightarrow
(\otimes\cdot(\mathcal{C}\times\otimes)\cdot a_{\mathcal{C}},
\varphi(\mathcal{C}\times\otimes)a_{\mathcal{C}})\\
&:& ((\mathcal{C}\times\mathcal{C})\times\mathcal{C},
(F\times\mathcal{C})\times \mathcal{C})\longrightarrow (\mathcal{C}, F)\\
r&:& (\otimes\cdot(\mathcal{C}\times\delta_{I})\cdot r^{-1}_{\mathcal{C}},
\separ\varphi(\mathcal{C}\times\delta_{I})\separ
r^{-1}_{\mathcal{C}})\longrightarrow (1_{\mathcal{C}}, 1_{\F}):(\mathcal{C},
F)\longrightarrow (\mathcal{C}, F)
\end{eqnarray*}

\item [3.-] Right actions on the Kleisli category, $\mathcal{C}_{\F}$,
  $\boxtimes:\mathcal{C}_{\F}\times\mathcal{C}\longrightarrow \mathcal{C}_{\F}$
  such that the following diagrams of morphisms and surfaces commute 

\begin{equation*}\Label{1405241922}
\xy<1cm,0cm>
\POS (0,20) *+{\mathcal{C}\times\mathcal{C}} = "a11",
\POS (20,20) *+{\mathcal{C}} = "a12",
\POS (0,0) *+{\mathcal{C}_{\F}\times\mathcal{C}} = "a21",
\POS (20,0) *+{\mathcal{C}_{\F}} = "a22",

\POS "a11" \ar^-{\otimes} "a12",
\POS "a12" \ar^{G_{\F}} "a22",
\POS "a11" \ar_{G_{\F}\times\mathcal{C}} "a21",
\POS "a21" \ar_-{\boxtimes} "a22",

\endxy
\end{equation*}

\begin{equation*}\Label{1405241923}
\begin{array}{ccc}
\xy<1cm,0cm>
\POS (0,15) *+{(a)}
\POS (0, 0) *+{\mathcal{C}^{2}\times\mathcal{C}} = "c1",
\POS (7, 0) *+{\phantom{(\mathcal{C})}}= "c3",
\POS (30, 0) *+{\mathcal{C}} = "c2", 
\POS (11, 4) = "a",
\POS (11, -4) = "b",
\POS (14,-36) = "c",
\POS (14,-44) = "d",

\POS (0,-40) *+{(\mathcal{C}_{\F}\times\mathcal{C})\times\mathcal{C}} = "d1",
\POS (30,-40) *+{\mathcal{C}_{\F}} = "d2",

\POS "c1" \ar@/^1.5pc/^{\otimes\cdot (\otimes\times \mathcal{C})} "c2",
\POS "c1" \ar@/_1.5pc/_{\otimes\cdot (\mathcal{C}\times\otimes)\cdot a_{\mathcal{C}}} "c2",
\POS "a" \ar^{\phantom{11}a} "b",

\POS "d1" \ar@/^1.5pc/^{\boxtimes\cdot (\boxtimes\times \mathcal{C})} "d2",
\POS "d1" \ar@/_1.5pc/_{\boxtimes\cdot (\mathcal{C}_{\F}\times\otimes)\cdot a_{\mathcal{C}_{\ast}}} "d2",
\POS "c" \ar^{\phantom{11}\tilde{a}} "d",

\POS "c1" \ar_{(G_{\F}\times\mathcal{C})\times\mathcal{C}} "d1",
\POS "c2" \ar^{G_{\F}} "d2",
\endxy & 
\xy<1cm,0cm>
\POS (0, 0) *+{\phantom{123}}
\endxy
&
\xy<1cm,0cm>
\POS (0,15) *+{(b)}
\POS (0, 0) *+{\mathcal{C}} = "c1",
\POS (7, 0) *+{\phantom{(\mathcal{C})}}= "c3",
\POS (30, 0) *+{\mathcal{C}} = "c2", 
\POS (11, 4) = "a",
\POS (11, -4) = "b",
\POS (11,-36) = "c",
\POS (11,-44) = "d",

\POS (0,-40) *+{\mathcal{C}_{\F}} = "d1",
\POS (30,-40) *+{\mathcal{C}_{\F}} = "d2",

\POS "c1" \ar@/^1.5pc/^{\otimes\cdot (\mathcal{C}\times\delta_{I})\cdot r^{-1}_{\mathcal{C}}} "c2",
\POS "c1" \ar@/_1.5pc/_{1_{\mathcal{C}}} "c2",
\POS "a" \ar^{\phantom{12}r} "b",

\POS "d1"
\ar@/^1.5pc/^{\boxtimes\cdot(\mathcal{C}_{\F}\times\delta_{I})\cdot r^{-1}_{\mathcal{C}_{\F}}} "d2",
\POS "d1" \ar@/_1.5pc/_{1_{\mathcal{C}_{\F}}} "d2",
\POS "c" \ar^{\phantom{12}\tilde{r}} "d",

\POS "c1" \ar_{G_{\F}} "d1",
\POS "c2" \ar^{G_{\F}} "d2",
\endxy
\end{array}
\end{equation*}

\end{enumerate}

\end{1505241921}

We left to the reader the writing of dual statements, \emph{i.e.} the ones that
corresponds to the Eilenberg-Moore category, where the direction of the
natural transformations are inverted, for example $\widehat{\sigma}^{r}_{A,B}:
F(A\otimes B)\longrightarrow A\otimes FB$.

\section{Functor Algebras}

Check Proposition II.1.1 in \cite{duej_kaee} and \cite{mewi_nobh} for this section. We define the category of $H$-left functor algebras for a given monad $(\mathcal{D}, H, \mu^{\Hh}, \eta^{\Hh})$. 

\newtheorem{1407022006}{Definition}[section]
\begin{1407022006}

The category of left $H$-functor algebras, for the pair $(\mathcal{C},\mathcal{D})$, denoted as ${}_{\Hh}\mathcal{F}$ or ${}_{\Hh}\mathcal{M}$ is defined as follows. The objects are given by $(J,\lambda_{\J})$, where $J:\mathcal{C}\longrightarrow\mathcal{D}$ is a functor and $\lambda_{\J}: HJ\longrightarrow J$ is a natural transformation such that the following diagrams commute 

\begin{equation}\Label{1407022029}
\begin{array}{cc}
\xy<1cm,0cm>
\POS (0,20) *+{HHJ} = "a11",
\POS (27,20) *+{HJ} = "a12",
\POS (0,0) *+{HJ} = "a21",
\POS (27,0) *+{J} = "a22",

\POS "a11" \arrow^{\mu^{\Hh}J} "a12",
\POS "a12" \arrow^{\lambda_{\J}} "a22",
\POS "a11" \arrow_{H\lambda_{\J}} "a21",
\POS "a21" \arrow_{\lambda_{\J}} "a22", 

\endxy &

\xy<1cm,0cm>
\POS (0,20) *+{J} = "a11",
\POS (20,20) *+{HJ} = "a12",
\POS (20,0) *+{J} = "a22",

\POS "a11" \arrow^{\eta^{\Hh}J} "a12",
\POS "a12" \arrow^{\lambda_{\J}} "a22",
\POS "a11" \arrow_{1_{\J}} "a22", 

\endxy
\end{array}
\end{equation}

A \emph{morphism of functor algebras}, is a natural transformation, $\theta:(J,\lambda_{\J})\longrightarrow (K,\lambda_{\K})$,
$\theta:J\longrightarrow K$ such that the following diagram commute

\begin{equation}\Label{1407022032}
\xy<1cm,0cm>
\POS (0,20) *+{HJ} = "a11",
\POS (25,20) *+{HK} = "a12",
\POS (0,0) *+{J} = "a21",
\POS (25,0) *+{K} = "a22",

\POS "a11" \arrow^{H\theta} "a12",
\POS "a12" \arrow^{\lambda_{\K}} "a22",
\POS "a11" \arrow_{\lambda_{\J}} "a21",
\POS "a21" \arrow_{\theta} "a22",

\endxy
\end{equation}

\end{1407022006}

We realize that the diagrams given by \eqref{1407022029}, for a left
$H$-functor algebra, account for a monad morphism of the form
$(J,\lambda_{\J}):(\mathcal{C}, 1_{\mathcal{C}})\longrightarrow (\mathcal{D},
H)$. In the same way, the commutative diagram for a  morphism of left
$H$-functor algebras, as in \eqref{1407022032},  account for a monad
transformation $\theta:(J,\lambda_{\J})\longrightarrow(K,
\lambda_{\K}):(\mathcal{C}, 1_{\mathcal{C}})\longrightarrow(\mathcal{D}, H)$.\\

Using the isomorphism for the Eilenberg-Moore 2-adjunction, given by \eqref{1403191902}, the category ${}_{\Hh}\mathcal{F}$ is isomorphic to the following category, named possibly as \emph{category of liftings to} $\mathcal{D}^{\Hh}$, for the pair $(\mathcal{C},\mathcal{D})$. The objects of such category are functor pairs $(J, \hat{J}\separ)$ such that they complete to an adjunction morphism, in $\Adj_{\R}({}_{2}Cat)$, of the form $(J, \hat{J}\separ): 1_{\mathcal{C}}\dashv 1_{\mathcal{C}}\longrightarrow D^{\Hh}\dashv U^{\Hh}$. That is to say, the following diagram commutes

\begin{equation*}
\xy<1cm,0cm>
\POS (0, 0) *+{\mathcal{C}} = "c1",
\POS (20, 0) *+{\mathcal{D}} = "c2", 
\POS (0, -20) *+{\mathcal{C}^{1_{\mathcal{C}}}} = "d1",
\POS (20, -20) *+{\mathcal{D}^{\Hh}} = "d2", 
\POS "c1" \ar^{J} "c2",
\POS "d1" \ar^{1^{1_{\mathcal{C}}}} "c1",
\POS "d2" \ar_{U^{\Hh}} "c2",
\POS "d1" \ar_{\hat{J}} "d2",
\endxy 
\end{equation*}

\emph{i.e.}

\begin{equation*}
\xy<1cm,0cm>
\POS (0, 0) *+{\mathcal{C}} = "c1",
\POS (20, 0) *+{\mathcal{D}} = "c2", 
\POS (20, -18) *+{\mathcal{D}^{\Hh}} = "d2", 
\POS "c1" \ar^{J} "c2",
\POS "d2" \ar_{U^{\Hh}} "c2",
\POS "c1" \ar_{\hat{J}} "d2",
\endxy 
\end{equation*}

The morphisms of such category are the usual morphisms of adjunctions
$(\alpha, \beta): (J, \hat{J}\separ)\longrightarrow
(K,\hat{K}):1_{\mathcal{C}}\dashv 1_{\mathcal{C}}\longrightarrow D^{\Hh}\dashv
U^{\Hh}$. We then proved the following theorem.

\newtheorem{1406031326}{Theorem}[section]
\begin{1406031326}

There exists an isomorphism, natural on $\mathcal{C}$ and $(\mathcal{D}, H)$, between the following categories

\begin{enumerate}

\item [1.-] The category of left $H$-functor algebras ${}_{\Hh}\mathcal{F}$. 

\item [2.-] The category of liftings to $\mathcal{D}^{\Hh}$, for the pair $(\mathcal{C},\mathcal{D})$.

\end{enumerate}

\end{1406031326}

\begin{flushright}
$\square$
\end{flushright}

Dually, we have the category of right $H$-functor algebras, for the monad $(\mathcal{D}, H, \mu^{\Hh}, \eta^{\Hh})$, denoted as $\mathcal{F}_{H}$ or $\mathcal{M}_{\Hh}$. The objects are pairs $(J, \rho_{\J})$, where the natural transformation $\rho_{\J}:JH\longrightarrow J$ is such that it fulfills diagrams dual to those in \eqref{1407022029}. In the same (dual) way as before, this category is the same as the category $Hom_{\Mnd^{\scriptscriptstyle{\bullet}}({}_{2}Cat)}((\mathcal{D}, H), (\mathcal{C}, 1_{\mathcal{C}}))$. Therefore using the isomorphism \eqref{1311042204}, the previous category is isomorphic to the category named as \emph{extensions from} $\mathcal{D}_{\Hh}$, for the pair $(\mathcal{D},\mathcal{C})$. The objects of this category are pairs of functors $(J, \tilde{J}\separ)$ such that thery complete to an adjunction morphism $(J,\tilde{J}\separ):G_{\Hh}\dashv V_{\Hh}\longrightarrow 1_{\mathcal{C}}\dashv 1_{\mathcal{C}}$ in $\Adj_{\Ll}({}_{2}Cat)$. In particular, the following diagram commutes 

\begin{equation*}
\xy<1cm,0cm>
\POS (0, 0) *+{\mathcal{D}} = "c1",
\POS (20, 0) *+{\mathcal{C}} = "c2", 
\POS (0, -18) *+{\mathcal{D}_{\Hh}} = "d2", 
\POS "c1" \ar^{J} "c2",
\POS "d2" \ar_{\tilde{J}} "c2",
\POS "c1" \ar_{G_{\Hh}} "d2",
\endxy 
\end{equation*}

\noindent We also proved the following theorem 

\newtheorem{1406032115}[1406031326]{Theorem}
\begin{1406032115}

There exists an isomorphism, natural on $(\mathcal{D}, H)$ and $\mathcal{C}$, between the following categories

\begin{enumerate}

\item [1.-] The category of right $H$-functor algebras $\mathcal{F}_{\Hh}$. 

\item [2.-] The category of extensions from $\mathcal{D}_{\Hh}$, for the pair $(\mathcal{D},\mathcal{C})$.

\end{enumerate}

\end{1406032115}

\begin{flushright}
$\square$
\end{flushright}

\section{Conclusions and Future Work}

This survey has the objective to show how several situations for the theory of
monads are connected in a very simple way, through a 2-adjunction. Any person who has
tought a course on monads would agree that this structure, of a 2-adjunction,
can be used as an educational purpose in the sense that a simple structure can
account for several situations and which can spare the, otherwise cumbersome,
details of the proofs. \\

For future work, we have a few recommendations. The reader may find interesting
to extent the part of strong monads and actions over the Kleisli categories to
strong symmetrical monads and use the actions for the Eilenberg-Moore case. It
would be interesting also to contextualize the case of the monoidal liftings
and monoidal extensions according to the formal theory of monoidal monads, and the
\emph{standar} objects, given in \cite{zama_fotm}.\\

The reader may want to find more situations in the monad theory that can use the isomorphism provided by this pair of 2-adjunctions, the authors will certainly pursue this issue.

\begin{center}
\section*{Acknowledgement}
\end{center}

The third author would like to thank to the Consejo de Ciencia y Tecnolog\'ia (CONACYT) for partial  
financial support through the grant SNI-59154. We all would like to thank to the referee who made some interesting comments on an earlier version.

\end{document}